\documentclass{article}
\usepackage{amsmath}
\usepackage{amsthm}
\usepackage{tikz}
\usepackage[all]{xy}
\usepackage{amssymb}
\usepackage{mathtools}
\usepackage{mathabx}
\usepackage{algorithm}
\usepackage{algpseudocode}
\usepackage[toc]{appendix}
\usepackage{url}
\usepackage[nameinlink]{cleveref}

\DeclarePairedDelimiter\floor{\lfloor}{\rfloor}
\DeclareMathAlphabet{\mathbbold}{U}{bbold}{m}{n}
\usepackage[numbers,sort&compress]{natbib}
\usepackage{tikz-cd}
\usepackage[all]{xy}
\newtheorem{theorem}{Theorem}[section]
\newtheorem{lemma}[theorem]{Lemma}

\newtheorem{remark}[theorem]{Remark}
\newtheorem{corollary}[theorem]{Corollary}
\newtheorem{proposition}[theorem]{Proposition}

\theoremstyle{definition}
\newtheorem{definition}[theorem]{Definition}
\newtheorem{question}[theorem]{Question}
\newtheorem{example}[theorem]{Example}

\theoremstyle{remark}

\usepackage{graphicx} 
\newcommand{\ZZ}{\mathbf{Z}}
\newcommand{\QQ}{\mathbf{Q}}
\newcommand{\RR}{\mathbf{R}}
\newcommand{\FF}{\mathbf{F}}
\newcommand{\CC}{\mathbf{C}}

\newcommand{\un}{\underline}
\newcommand{\Bin}{\mathrm{Bin}}
\newcommand{\hBin}{\text{-}\mathrm{Bin}}
\newcommand{\Lie}{\mathrm{Lie}}
\newcommand{\Hom}{\mathrm{Hom}}
\newcommand{\val}{\mathrm{val}}

\newcommand{\an}{\mathrm{an}}

\newcommand{\Gal}{\mathrm{Gal}}

\newcommand{\Mod}{\mathrm{Mod}}

\newcommand{\la}{\mathrm{la}}

\newcommand{\A}{\mathrm{A}}

\newcommand{\Id}{\mathrm{Id}}
\newcommand{\han}{\text{-an}}
\newcommand{\bs}{\blacksquare}

\newcommand{\LS}{\mt{LS}}
\newcommand{\Koh}{\mt{Koh}}
\newcommand{\Spa}{\mt{Spa}}
\newcommand{\Solid}{\mt{Solid}}
\newcommand{\rig}{\mt{rig}}
\newcommand{\Cont}{\mt{Cont}}
\newcommand{\RHom}{\mt{RHom}}
\newcommand{\Bbs}{B_{\bs}}
\newcommand{\Bplusbs}{B^{+}_{\bs}}

\newcommand{\wt}{\widetilde}
\newcommand{\mc}{\mathcal}
\newcommand{\mt}{\mathrm}
\newcommand{\slope}{\mt{slope}}
\newcommand{\mfrk}{\mathfrak}
\newcommand{\bb}{\mathbb}
\newcommand{\wh}{\widehat}
\newcommand{\ul}{\underline}
\newcommand{\Fil}{\mt{Fil}}
\newcommand{\gr}{\mt{gr}}

\newcommand{\angles}[1]{\langle #1 \rangle}
\newcommand{\ra}{\rightarrow}
\newcommand{\Dh}{\mc{D}_{h \han}(G,B)}
\newcommand{\Dhplus}{\mc{D}_{h \han}(G,B^+)}
\newcommand{\Dhpluszero}{\mc{D}_{h \han}(G_0,B^+)}
\newcommand{\Dhzero}{\mc{D}_{h \han}(G_0,B)}
\title{Solid locally analytic representations in mixed characteristic}
\author{Gal Porat
\\{galporat1@gmail.com}}
\begin{document}

\maketitle

\begin{abstract}
The theory of locally analytic representations of $p$-adic Lie groups with $\QQ_p$-coefficients is a powerful tool in $p$-adic Hodge theory and in the $p$-adic Langlands program. This perspective reveals important differential structures, such as the Sen and Casimir operators.

 Rodríguez Camargo and Rodrigues Jacinto  developed in \cite{RJRC22} a solid version of this theory using the language of condensed mathematics. This provides more robust homological tools (comparison theorems, spectral sequences...) for studying these representations.

In this article, we extend the solid theory of locally analytic representations to a much broader class of mixed characteristic coefficients, such as $\FF_p((X))$ or $\ZZ_p[[X]]\angles{p/X}[1/X]$, as well as to semilinear representations. In the introduction, we explain how these ideas could relate to mixed characteristic phenomena in $p$-adic Hodge theory, extend eigenvarieties, and the Langlands program.
\end{abstract}

\maketitle

\tableofcontents

\section{Introduction}

This article is concerned with the locally analytic representation theory of $p$-adic Lie groups with mixed characteristic coefficients. Before stating the main results in $\mathsection$\ref{subsec:the main results}, we explain how these objects arise naturally and why it is interesting to study them. 

\subsection{Motivation from an example}
\label{subsec:motivation_from_an_example}
In the following example, we show how mixed characteristic locally analytic representations naturally come up in a very simple setting. Let $p$ be a prime and let $${\bf{D}}_{\QQ_p}=\Spa\ZZ_p[[T]]\times_{\ZZ_p}\QQ_p$$ be the open unit disc over $\QQ_p$. Consider the character $\lambda_T:x \mapsto (1+T)^x$ which takes values in $\mc{O}({\bf{D}}_{\QQ_p})$. This character is locally analytic, in the sense of $p$-adic representation theory, when its values are restricted to affinoid subspaces of ${\bf{D}}_{\QQ_p}$. For example, for every $m\in \ZZ_{\ge0}$, the closed disc $${\bf{D}}_{[0,p^{-m}]}=\{t:0\le |t|\le p^{-m}\}=\Spa\QQ_p\angles{T/p^m}$$
is an affinoid, and the character $\lambda_T:x\mapsto(1+T)^x$ from $\ZZ_p$ to $\QQ_p\angles{T/p^m}=\mc{O}({\bf{D}}_{[0,p^{-m}]})$ has a Taylor series expansion
\begin{align}
\label{eq:lambda_t_taylor_expansion_family}
\lambda_T(x) = \sum_{n=0}^\infty \frac{\log(1+T)^n}{n!}x^n.
\end{align}
The Banach norm of the coefficients $a_n=\log(1+T)^n/n!$ decays exponentially, which means $\lambda_T$ is a locally analytic representation of $\ZZ_p$. As a continuous rank $\ZZ_p$-representation, it is a module over the algebra of $\ZZ_p$-measures $\ZZ_p[[\ZZ_p]][1/p]$, and it being locally analytic implies this module structure extends to the larger algebra $\mc{D}^\la(\ZZ_p,\QQ_p)$ of locally analytic distributions. Furthermore, it implies the existence of  an action of the Lie algebra $\Lie(\ZZ_p)$. Since $d/dx((1+T)^x)|_{x=0}=\log(1+T),$ this essentially amounts to the "multiplication by $\log(1+T)$" operator.

The character $\lambda_T$ also makes sense after taking its mod $p$ reduction. It therefore makes sense to ask if one can extend the family $\lambda_T$ to include characteristic $p$ points as locally analytic characters as well. The answer is positive - let us try to make sense of this. The open unit disc ${\bf{D}}_{\QQ_p}$ is contained in the extended disc
$${\bf{D}}=\Spa(\ZZ_p[[T]],\ZZ_p[[T]])^\an,$$ in which we have included the analytic mod $p$ fiber.
This larger space includes pseudorigid subaffinoids as well, such as the pseudorigid discs
$${\bf{D}}_{[p^{-m},1]}=\{t:p^{-m}\le|t|\le1 \}=\Spa\ZZ_p([[T]]\angles{p/T^m}[1/T]),$$
in which the inequalities on $|t|$ are imposed in the opposite direction. 
The character $\lambda_T$ extends $\bf{D}$, and in particular to the pseudorigid discs ${\bf{D}}_{[p^{-m},1]}$. However, in the Banach ring $\ZZ_p([[T]]\angles{p/T^m}[1/T])$, the element $p$ is not invertible! It is a mixed characteristic Banach ring which is not a $\QQ_p$-algbera. Thus, $\lambda_T:\ZZ_p\to \mc{O}({\bf{D}}_{[p^{-m},1]})$ certainly does not have a Taylor series expansion
as in (\ref{eq:lambda_t_taylor_expansion_family}). Nevertheless, it makes sense to say that $\lambda_T$ is locally analytic, even on these pseudorigid discs, if we rephrase everything in terms of binomial expansions. The idea is to replace (\ref{eq:lambda_t_taylor_expansion_family}) with
\begin{align}
\label{eq:lambda_t_binomial_expansion}
\lambda_T(x)=\sum_{n=0}T^n\binom{x}{n}.
\end{align}
According to a theorem of Amice, a function valued in a $\QQ_p$-Banach algebra is locally analytic if and only if in its binomial expansion $\sum_{n=0}a_n\binom{x}{n}$, the coefficients $a_n$ tend to 0 exponentially. We use this as a \textit{definition} of a locally analytic function in the pseudorigid case. On $\mc{O}({\bf{D}}_{[0,p^{-m}]})$, where the topology is $p$-adic, this condition holds because of the equation $T = p^m\cdot (T/p^m)$. On the other hand, the topology of $\mc{O}({\bf{D}}_{[p^{-m},1]})$ is $T$-adic, so the equation (\ref{eq:lambda_t_binomial_expansion}) shows that $\lambda_T$ is locally analytic even on ${\bf{D}}_{[p^{-m},1]}$. One can then define an appropriate algebra of distributions $\mc{D}^\la(\ZZ_p,\ZZ_p[[T]]\angles{p/T}[1/T])$-coefficients, and $\lambda_T$ has a module structure over it, similarly to the case of $\QQ_p$-coefficients. We can even reduce mod $p$ and get a locally analytic character $\overline{\lambda}_T:\ZZ_p\to \FF_p((T))$. 


\subsection{Motivation from $p$-adic Hodge theory}
 In $p$-adic Hodge theory, one often encounters a module $M$ over a ring $R$ with a semilinear action of some $p$-adic Lie group $G$. There is a $G$-stable subring $R_0 \subset R$ with $R = \widehat{R_0}$, and one wants to know if $M$ can be descended to a $G$-stable module $M_0$ over $R_0$. We loosely call such a descent $M_0$ a decompletion of $M$.
 
 For example, Fontaine's theory attaches to every 
 $p$-adic representation $\rho$ of $\Gal(\overline{\QQ}_p/\QQ_p)$ with $\ZZ_p$-coefficients an étale $(\varphi,\ZZ_p^\times)$-module. These objects are modules over a certain power series ring ${\bf{A}}_{\QQ_p}$, equipped with a semilinear $(\varphi,\ZZ_p^\times)$-action. Using the decompletion theory of Cherbonnier-Colmez (\cite{cherbonnier1998representations}), one knows each such  $(\varphi,\ZZ_p^\times)$-module descends to a module over the smaller overconvergent subring ${\bf{A}}_{\QQ_p}^\dagger \subset {\bf{A}}_{\QQ_p}$ (here the $p$-adic Lie group is $\ZZ_p^\times$). Unlike the complete ring ${\bf{A}}_{\QQ_p}$, the ring ${\bf{A}}_{\QQ_p}^\dagger$ is contained in the Robba ring $\mc{R}_{\QQ_p}$ from the theory of $p$-adic differential equations, and by a theorem of Kedlaya we get from this a functor $\rho\mapsto D_{\mc{R}_{\QQ_p}}(\rho)$, associating to $\rho$ a  $(\varphi,\ZZ_p^\times)$-module over $\mc{R}_{\QQ_p}$ (\cite{kedlaya2004p}). Putting all of this together gives a beautiful link between the theory of $p$-adic Galois representations and $p$-adic analysis. This connection was used in the resolution of Fontaine's $p$-adic monodromy conjecture (\cite{berger2002representations}).

 In recent years, it has become clear there is an intimate relation between decompletion and the theory of locally analytic representations of $p$-adic Lie groups. Namely, it was recognized that that decompletion is usually of the form
$$M =R \otimes_{R^\la}M^\la,$$
and hence can be proved and studied with techniques of locally analytic representation theory. 
 This idea was first noticed in the paper of Berger and Colmez (\cite{berger2016theorie}) which studied Sen theory via the functor of locally analytic vectors. It was applied further in various settings where $R$ is a $\QQ_p$-Banach algebra (see \cite{berger2016multivariable}, \cite{gao2021locally}, \cite{porat2024locally}, \cite{pan2022locally}, \cite{camargo2022geometric}, and others). 
 However, decompletion in $p$-adic Hodge theory already occurs on an integral level: indeed, the ring ${\bf{A}}_{\QQ_p}$ mentioned above is not a $\QQ_p$-algebra. A first step in understanding decompletion via locally analytic methods in an integral setting was taken by Berger and Rozenstazjn. They showed in a mod $p$ setting that the field of norms of a perfectoid extension can be understood using locally analytic vectors (\cite{BR22a,BR22b}). The author was inspired by this work in \cite{Po24} to study decompletion using locally analytic representations with mixed characteristic coefficients. This showed how to interpret the overconvergent ring ${\bf{A}}_{\QQ_p}^\dagger$, where $p$ is not invertible, in terms of locally analytic vectors. The methods and definitions of our article were a little ad hoc, and their shortcomings motivated a more streamlined development of the theory.

\subsection{Motivation from the theory of extended eigenvarieties}

Recall that an eigenvariety is a certain rigid analytic space $\mc{E}^\rig$ which parametrizes overconvergent $p$-adic forms of finite slope for a given reductive group $G$ over $\QQ$. A feature of this eigenvariety is a map from the eigenvariety to the space of weights $\mc{W}^\rig$, which encodes the $p$-adic weight (or weights) of the $p$-adic form. The first example of an eigenvariety is the eigencurve $\mc{C}^\rig$ constructed by Coleman and Mazur in the setting of $G=\mt{GL}_2$ (\cite{coleman1998eigencurve}). In this case, $\mc{W}^{\rig} = \Spa(\ZZ_p[[\ZZ_p^\times]])^\rig\cong \coprod_{(\ZZ/p)^\times}{\bf{D}}_{\QQ_p}$. The question of constructing eigenvarieties is an important problem which has a rich history. It can be done by at least two methods: the method of overconvergent cohomology (\cite{hansen2017universal}) and the locally analytic Jacquet functor method of Emerton (\cite{emerton2006interpolation}).

In the beautiful article "Le Halo Spectral" (\cite{andreatta2018halo}) by Andreatta, Iovita and Pilloni, the authors construct spaces of mixed-characteristic overconvergent modular forms, and as a consequence, an extended eigencurve $\mc{C}\supset \mc{C}^\rig$ with a map to an extended weight space $\mc{W}\supset\mc{W}^\rig$.  This weight space $\mc{W}$ is none other than $\coprod_{(\ZZ/p)^\times}{\bf{D}}$ for the extended disc ${\bf{D}}$ appearing in $\mathsection$\ref{subsec:motivation_from_an_example}. Later, wide generalizations of the construction of (\cite{andreatta2018halo}) were obtained by Johansson-Newton (\cite{johansson2019extended}) and Gulotta (\cite{gulotta2019equidimensional}), both by adapting the method of overconvergent cohomology to mixed characteristic.

With that said, following natural question is due to Rebecca Bellovin. A positive answer would establish a missing link between the theory of extended eigenvarieties and locally analytic representation theory in mixed characteristic, which already exists for $\QQ_p$-coefficients.

\begin{question}
Is it possible to define a mixed-characteristic version of Emerton's
locally analytic Jacquet functor, in order to construct extended
eigenvarieties from completed cohomology?
\end{question}

Recall the idea in the $\QQ_p$ setting for $G=\mt{GL}_2$: for some fixed tame level $K^p$, one takes the first completed cohomology of modular curves $$\tilde{\mt{H}}^1(K^p,\ZZ_p)=\varprojlim_n\varinjlim_{K_p\subset \mt{GL}_2(\QQ_p)}\mt{H}^1(Y_{K^pK_p}(\bb{C}),\ZZ/p^n).$$
This space is $p$-adically complete and has an action of $\mt{GL}_2(\QQ_p)$. A key object in Emerton's representation-theoretic construction of the eigencurve is the locally analytic Jacquet module of $$(\tilde{\mt{H}}^1(K^p,\ZZ_p)\otimes_{\ZZ_p}\QQ_p)^\la.$$
Recall that $\tilde{\mt{H}}^1(K^p,\ZZ_p)$ is a module over the big Hecke algebra $\bb{T}(K^p)$, hence over the ring of weights $\ZZ_p[[\ZZ_p^\times]]\cong \ZZ_p[\ZZ/p^\times][[T]]$. Perhaps considering a locally analytic Jacquet module of 
$$(\tilde{\mt{H}}^1(K^p,\ZZ_p)\otimes_{\ZZ_p[[T]]}\ZZ_p[[T]]\angles{p/T}[1/T])^\la$$
instead could lead to a construction of the extended eigencurve.

\subsection{Motivation from the Langlands program}

The Langlands program connects automorphic representations and Galois representations. Recall the following rich picture we have in the setting of $\mt{GL}_2$ with $L$-coefficients, for $L$ a finite extension of $\QQ_p$ (proven under some mild conditions):

1. Every overconvergent eigenform $f$ of finite slope with $L$-coefficients gives rise to a Galois representation $\rho:\Gal(\overline{\QQ}/\QQ)\to \mt{GL}_2(L)$ which is odd, 
    unramified almost everywhere with $D_{\mc{R}_{L}}(\rho|_{\Gal(\overline{\QQ}_p/\QQ_p)})$ trianguline (\cite[Theorem 6.3]{kisin2003overconvergent}). Every Galois representation satisying these conditions arises in this way (\cite[Corollary 1.2.2]{emerton2011local}).

2. There exists a correspondence between admissible unitary continuous Banach representations of $\mt{GL}_2(\QQ_p)$ with $L$ coefficients and Galois representations $\Gal(\overline{\QQ}_p/\QQ_p)\to \mt{GL}_2(L)$ (\cite{colmez2010representations}).

3. There exists a correspondence between certain admissible locally analytic representations of $\mt{GL}_2(\QQ_p)$ with $L$-coefficients and rank 2 $(\varphi,\ZZ_p^\times)$-modules over $\mc{R}_L$ (\cite{colmez2016representations}), compatible with the previous correspondence. Trianguline étale $(\varphi,\ZZ_p^\times)$-modules correspond to a representation with non vanishing Jacquet module (\cite[Théorème 1.1.]{dospinescu2011equations}).

In 2 and 3 above the correspondences satisfy certain nice properties we do not wish to spell out here.


\begin{question}
Do we have a similar picture in the case where $L$ is a finite extension of $\FF_p((T))$?
\end{question}

At least for part 1 above the answer is positive under some mild conditions by results of Bellovin (\cite{bellovin2024cohomology,bellovin2024modularity}). In that case, we see the overconvergent eigenforms in characteristic $p$ appearing in the boundary of compactification the space of characteristic 0 overconvergent eigenforms. It is natural to ask the same, in a local context, for the spaces of representations appearing  in the analytic categorical $p$-adic Langlands correspondence (\cite[$\mathsection{6.2}$]{emerton2022introduction}).

\begin{question}
1. Is there a compactification of the analytic stack of $(\varphi,\Gamma)$-modules (resp. trianguline $(\varphi,\Gamma)$-modules) $\mathfrak{X}_{\mt{GL}_d}$ (resp.  $\mathfrak{X}_{\mt{GL}_d,\mt{tri}}$) which incorporates mixed characteristic and characteristic $p$ points?

2. If the answer to 1 is positive - is there an extension of the conjectural categorical $p$-adic Langlands correspondence to this compactification?
\end{question}

\subsection{The main results}
\label{subsec:the main results}
The theory of $p$-adic locally analytic representations was initiated in a series of papers of Schneider and Teitelbaum (\cite{schneider_teitelbaum2002algebras,ST01,ST02locally,ST02}) and developed by many people. Their foundations were also reworked by Emerton (\cite{emerton2017locally}) and by Rodrigues Jacinto and Rodríguez Camargo in the solid context (\cite{RJRC22,RJRC23}). 

In this article we build the basic theory of solid locally analytic representations in mixed characteristic, extending results in the $p$-adic case due to \cite{RJRC22}. Let $G$ be a compact $p$-adic Lie group, and let $B$ be a mixed characteristic Banach ring (see Definition \ref{def:banach_pair}), satisfying some conditions which capture most examples of interest.\footnote{We require that $B$ comes from a Banach pair $(B,B^+)$ which has slope $\le 1$ and is residually of finite type, see Definitions \ref{def:slope} and  \ref{def:residually_of_finite_type}.} Examples of Banach rings satisfying\footnote{Nonexamples are $\CC_p, \FF_p((T^{1/p^\infty}))$ or $\ZZ_p[[X]]\angles{p^2/X}[1/X]$.} the assumptions are $\QQ_p,\FF_p((T))$ or $\ZZ_p[[X]]\angles{p/X}[1/X]$. A Banach ring $B$ give rise to an analytic ring $\Bbs$ in the sense of condensed mathematics, see $\mathsection{\ref{subsec:banach_pairs_as_an_rings}}$.

The theory works as follows. We have the Iwasawa $B$-algebra $$\Bbs[G]=\varprojlim_N\Bbs^+[G/N].$$ The category of solid $\Bbs[G]$-modules is the same as the category of solid $B$-modules with a continuous $G$-action. In particular, it contains the category of continuous $G$-representations on Banach spaces over $B$. Using binomial expansions, we define rings of analytic functions $\mc{C}_{h \han}(G,B)$ and analytic distributions $\mc{D}_{h \han}(G,B)$ containing $\Bbs[G]$. For a Banach space $V$ over $B$ we can define its $h$-analytic vectors to be these elements with locally analytic $G$-action, i.e.
$$V_{h \han}=\mt{RHom}_{\Bbs[G]}(B,\mc{C}_{h \han}(G,B)\otimes_{\Bbs} V).$$
This definition can be extended to general solid $\Bbs[G]$-modules and the derived category $D(\Bbs[G])$. The $h$-analytic vectors are a module over the distribution algebra $\mc{D}_{h \han}(G,B)$. A representation $V$ is said to be $h$-analytic if the map $V_{h \han}\to V$ is an isomorphism.

Our first main theorem is the following.

\begin{theorem}[Theorem \ref{thm:full_faithfullness}]
\label{thm:full_faithfulness_intro}
1. The category of solid $h$-analytic representations is a full subcategory of the solid $\Bbs[G]$-modules.

2. The category of $h$-analytic $G$-representations over $B$ is equal to the category of solid modules over $\mc{D}_{h \han}(G,B)$.

3. A complex $C \in D(B_{\bs}[G])$ is  $h$-analytic  if and only if for all $n \in \ZZ$ the cohomologies $H^n(C)$ are $h$-analytic. Equivalently, $C$ is in the essential image of $D(\mc{D}_{h \han}(G,B)) \to D(\Bbs[G])$.
\end{theorem}

Given a complex $C \in D(\Bbs[G])$ we can define the derived locally analytic vectors $$C^{\la} = \varinjlim_h C^{h \han}.$$ Our next theorem is a comparison between the $G$-cohomology of $C$ and of its locally analytic vectors. 

\begin{theorem}[Theorem \ref{thm:cohomological_comparison}]
\label{thm:cohomological_comparision_intro}
We have
$$\mt{RHom}_{\Bbs[G]}(B,C) = \mt{RHom}_{\Bbs[G]}(B,C^\la).$$
\end{theorem}

As a consequence, we prove a (generalization of) conjecture 3.4 of \cite{Po24}. 

\begin{corollary}
\label{cor:spectral_seq_intro}
Let $V$ be a solid $\Bbs[G]$-module and let $\mt{R}^i_{\la}(V)$ be its derived locally analytic vectors. Then for $i \ge 0$ there exists a spectral sequence 
$$E_2^{i,j} = \mt{Ext}_{B_{\bs}[G]}^i(B,\mt{R}^j_{\la}(V)) \Longrightarrow \mt{Ext}_{B_{\bs}[G]}^{i+j}(B,V).$$
\end{corollary}

In particular, if there are no higher locally analytic vectors, $V^\la$ can be used to compute continuous cohomology:
\begin{corollary}
\label{cor:loc_ana_vecs_compute_cohomology}
If $\mt{R}^i_{\la}(V) = 0$ for $i \ge 1$, then for $i \ge 0$
$$\mt{H}^i(G,V) = \mt{H}^i(G,V^{\la}).$$
\end{corollary}

We also have a version of Lazard's theorem which compares between continuous and analytic cohomology.

\begin{theorem}[Theorem \ref{thm:lazard_comparison}]
\label{thm:lazard_intro}
Let $C\in D(\Bbs[G])$ be $h$-analytic. Then 
$$\mt{RHom}_{\Bbs[G]}(B,C)=\mt{RHom}_{\mc{D}_{h \han}(G,B)}(B,C).$$
\end{theorem}

In practice, one is sometimes interested in cases where $G$ acts on $B$ nontrivially, in other words one wants to consider semilinear $G$-representations with $B$-coefficients. For example, this is the case for $(\varphi,\ZZ_p^\times)$-modules mod $p$, where one has $G=\ZZ_p^\times$ and $B=\FF_p((X))$, with the action $a(X)=(1+X)^a-1.$ To deal with this case also, we define twisted rings $\Bbs[G]'$ and $\mc{D}_{h \han}(G,B)'$ where $G$ acts on $B$, i.e. we have the identity $[g]\cdot b=g(b)\cdot[g]$. We say that the action of $G$ on $B$ is locally analytic if\footnote{See Remark \ref{rem:unambiguity_of_analytic_action} for this terminology.} $(g-1)(\varpi^nB^+)\subset \varpi^{n+1}B^+$ for $n \in \ZZ_{\ge 0}$, for all $g$ in some open subgroup of $G$ and for $B^+\subset B$ the open unit ball. Note that this condition generalizes  the linear case, where $G$ acts trivially on $B$. Throughout, we actually work in this more general semilinear setting, and so we get the following extension.

\begin{theorem}
Suppose that the action of $G$ on $B$ is locally analytic. Then Theorem \ref{thm:full_faithfulness_intro}, Theorem \ref{thm:cohomological_comparision_intro}, Corollary \ref{cor:spectral_seq_intro}, Corollary \ref{cor:loc_ana_vecs_compute_cohomology} and Theorem \ref{thm:lazard_intro} hold for modules over the twisted rings $\Bbs[G]'$ and $\mc{D}_{h \han}(G,B)'$.
\end{theorem}

 \begin{remark} This semilinearity is what causes most of the technical problems we have to deal with in this paper, because $B$ is no longer central in the nasty rings $\Bbs[G]'$ and $\mc{D}_{h \han}(G,B)'$. The trick is to notice that, for $G$ uniform, and for certain filtrations, we have
$\gr(\Bbs[G]')=\gr(\Bbs[G])$ and $\gr(\mc{D}_{h \han}(G,B)')=\gr(\mc{D}_{h \han}(G,B))$. By using graded techniques going back to Schneider and Teitelbaum, structuring the arguments in the right way allows us to reduce the heavy lifting from the semilinear case to the linear case. 
 \end{remark}

\subsection{Further directions}

Due to time limitations, in this paper we assume $G$ is compact and we do not deal with admissible representations and smooth representations. It should be possible to extended the theory to include these. Removing the condition $\slope(B,B^+) \le 1$ is also desirable as ideally one would want to have a theory for pseudorigid coefficients such as $B=\ZZ_p[[X]]\angles{p^m/X^n}[1/X]$. Finally, it seems crucial to study whether there is some mysterious Lie algebra action in mixed characteristic. It is currently still missing from the picture. Indeed, even in the example discussed in $\mathsection{\ref{subsec:motivation_from_an_example}}$, the Lie algebra acts by multiplication by $\log(1+T)$ when $p\ne 0$. This series has an essential singularity at $p=0$, so the obvious attempt to extend the Lie algebra action does not work.

\subsection{Structure of the article}
In $\mathsection{2}$ we give some reminders on solid condensed mathematics. Next in $\mathsection{3}$ we set up the functional analysis over mixed characteristic Banach rings, and in $\mathsection{4}$ we give the definitions of the analytic functions and rings. In $\mathsection{5}$ we define the categories of continuous, analytic and locally analytic representations. Finally, in $\mathsection{6}$ we prove the main theorems.

\subsection{Notation and conventions}

They are the following.

1. If $\un{n} \in \ZZ^d$ we write $|\un{n}|$ for $\sum_{1\leq i \leq d}|n_i|$.

2. By a valuation on a ring $R$, we mean a map $\val_R:R \rightarrow (-\infty, \infty] $ satisfying
the following properties for $x, y \in R$:

(i). $\val_R(x) = \infty$ if and only if $x = 0$ (i.e. $R$ is separated);

(ii). $\val_R(xy) \geq \val_R(x) + \val_R(y)$;

(iii). $\val_R(x+y) \geq \min(\val_R(x), \val_R(y))$.

This definition can be extended in an obvious way
to an $R$-module $M$.

3. If $X$ is a topological space, we let $\ul{X}$ denote its associated condensed set. We sometimes identify $X$ and $\ul{X}$ when $X \mapsto \ul{X}$ is fully faithful, but we usually indicate before doing so.

\subsection{Acknowledgments}

I would like to thank Shay Ben-Moshe, Shai Keidar,  Juan Esteban Rodríguez Camargo, Joaquín Rodrigues Jacinto and Rustam Steingart for helpful comments throughout the writing of this project. It will be the obvious to the reader that this article leans heavily on the beautiful theory developed in \cite{RJRC22}, which was a pleasure to learn.

\section{Brief reminders on solid condensed mathematics}

In this section, we give a brief reminder of some definitions from the theory of solid condensed mathematics introduced by Dustin Clausen and Peter Scholze that are needed for this article. We assume the reader is familiar with the very basics of the theory, such as the definition of a condensed object. We refer the reader to \cite{CS19con}, \cite{CS19an}, \cite[$\mathsection{2}$]{Ma22} for a thorough treatment and \cite[$\mathsection{2}$]{RJRC22}, \cite[$\mathsection{2}$]{andreychev2021pseudocoherent} for a summary. 

\subsection{Analytic rings}

We recall the definition of an analytic ring following \cite[Definition 2.3.1]{Ma22} and \cite[Definition 2.3.10]{Ma22}, see also \cite[Definition 7.1]{CS19con}, \cite[Definition 7.4]{CS19con}, \cite[Definition 6.12]{CS19an}.
\begin{definition}
1. An uncompleted pre-analytic ring $\mc{A}=(\ul{\mc{A}},\mc{M})$ is the data of condensed ring $\ul{\mc{A}}$ together with a functor $S\mapsto \mc{M}[S]$ from extremally disconnected set to condensed $\ul{\mc{A}}$-modules, taking finite disjoint unions to products, and a natural transformation $S \ra \mc{M}[S]$.

2. An uncompleted analytic ring is an uncompleted pre-analytic ring $(\ul{\mc{A}},\mc{M})$ such that for any complex  $$...\ra C_1 \ra C_0 \ra 0$$ of condensed $\ul{\mc{A}}$-modules, such that all $C_\bullet$ are direct sums of objects of the form
$\mc{M}[S_i]$ for varying extremally disconnected $S_i$, the map\footnote{Here $\ul{\mc{A}}[S]$ is the free condensed $\ul{\mc{A}}$-module on $S$.}
$$\ul{\mt{RHom}}(\mc{M}[S], C) \ra \ul{\mt{RHom}}(\ul{\mc{A}}[S], C)$$
of condensed abelian groups is an isomorphism for all extremally disconnected sets $S$.

3. An analytic ring is an uncompleted analytic ring $\mc{A}$ such that the map $\ul{\mc{A}}\to \mc{M}[*]$ is an isomorphism.
\end{definition}

\begin{remark}
\label{rem:rem_on_an_rings}
1. If $\mc{A}$ is an uncompleted analytic ring, the functor $S\mapsto \mc{M}[S]$ can be extended to profinite sets. A priori, $\mc{M}[S]$ only lies in the derived category, but by \cite[Proposition 2.11]{andreychev2021pseudocoherent}, it is actually static whenever $\mc{A}$ is an uncompleted analytic ring over $\ZZ_{\bs}$. This will always be the case for any uncompleted analytic ring we encounter. 

2. When $\mc{A}$ is an analytic ring, we denote $\mc{A}[S] := \mc{M}[S],$ removing $\mc{M}$ from the notation.

3. The above definition of an uncompleted analytic ring is the same as Clausen and Scholze's definition of an analytic ring, and the above definition of an analytic ring is the same as their definition of a normalized analytic ring (\cite[Definition 12.9]{CS19an}).
\end{remark}

\begin{example}
\label{ex:analytic_rings}
1. Let $(A,A^+)$ be a complete Huber pair. Theorem 3.28 of \cite{andreychev2021pseudocoherent} constructs an analytic ring\footnote{That $(A,A^+)_{\bs}$ is an analytic ring and not merely an uncompleted analytic ring is explained in \cite[Remark 13.17]{CS19an}.} $(A,A^+)_{\bs}$ with underlying condensed ring $\ul{A}$. By Proposition 3.34 of loc. cit., the association $(A,A^+)\mapsto (A,A^+)_{\bs}$ is fully faithful. When $A^+=A^\circ$ is the subring of powerbounded elements, we simply write $A_{\bs}$ for $(A,A^+)_{\bs}.$

2. Let $\mc{A}$ be an analytic ring such that $\ul{\mc{A}}$ is static. Let $G$ be a condensed group acting on $\mc{A}$, i.e., a condensed map $G\times \ul{\mc{A}}\ra \ul{\mc{A}}$ satisfying the usual compatibility conditions, such that for every $g \in G$, the action map $\mt{act}(g):\ul{\mc{A}}\ra\ul{\mc{A}}$ induces a map of analytic rings $\mt{act}(g):\mc{A}\ra\mc{A}$. Then there is a (in general non static) analytic ring $\mc{A}[G]'$ whose ring structure is defined by the twisting identities $g\cdot a = g(a) \cdot g$ (see \cite[Definition 3.4.1]{Ma22}). If the action of $G$ on $\mc{A}$ is trivial, we simply write $\mc{A}[G]$ for $\mc{A}[G]'$. 
\end{example}

\subsection{Solid modules}

\begin{definition}
Let $\mc{A}$ be an uncompleted analytic ring. 

1. An $\ul{\mc{A}}$-module $M$ is said to be a solid $\mc{A}$-module if for every extremally disconnected set, the map
$$\mt{Hom}_{\ul{\mc{A}}}(\mc{M}[S],M)\ra \mt{Hom}_{\ul{\mc{A}}}(\ul{\mc{A}}[S],M)$$
is an isomorphism. 

2. Similarly, a complex $C$ in the derived category of condensed $\ul{\mc{A}}$-modules is said to be solid over $\mc{A}$ 
 if for every extremally disconnected set, the map
$$\mt{RHom}_{\ul{\mc{A}}}(\mc{M}[S],C)\ra \mt{RHom}_{\ul{\mc{A}}}(\ul{\mc{A}}[S],C)$$
is an isomorphism.
\end{definition}

\begin{theorem}[{\cite[Proposition 7.5]{CS19con}}]
\label{thm:solid_modules}
Let $\mc{A}$ be an uncompleted analytic ring.

1. The category $\mt{Mod}^{\mt{solid}}_{\mc{A}}$ of solid $\mc{A}$-modules is a full subcategory of the category of condensed $\ul{\mc{A}}$-modules. It is stable under all limits, colimits and extensions. Objects of the form $\mc{M}[S]$ for $S$ extremally disconnected are a family of compact projective generators of this category. The inclusion functor admits a left adjoint "solidification" functor
$$\mt{Mod}_{\ul{\mc{A}}}\to \mt{Mod}^{\mt{solid}}_{\mc{A}}, M \mapsto M\otimes_{\ul{A}} \mc{A}$$
which is colimit preserving and maps $\ul{\mc{A}}[S]$ to $\mc{M}[S]$.
There is a unique monoidal tensor product $\otimes_{\mc{A}}$ making the functor $M \mapsto M\otimes_{\ul{\mc{A}}} \mc{A}$ monoidal.

2. The derived category $D(\mc{A})$ of $\mt{Mod}^{\mt{solid}}_{\mc{A}}$ is a full subcategory of the derived category $D(\ul{\mc{A}})$ of condensed $\ul{\mc{A}}$-modules. It consists of these complexes which are solid over $\mc{A}$. A complex is solid over $\mc{A}$ if and only if each of its cohomologies is a solid module over $\mc{A}$.  The inclusion functor admits a left adjoint "solidification" functor
$$D(\ul{\mc{A}})\to D(\ul{\mc{A}}), C \mapsto C\otimes_{\ul{\mc{A}}}^L \mc{A}$$
which is colimit preserving and is the left derived functor of $M\mapsto M\otimes_{\ul{\mc{A}}}\mc{A}.$ There is a unique monoidal tensor product $\otimes_{\mc{A}}^L$ making the functor $C \mapsto C\otimes_{\ul{\mc{A}}}^L \mc{A}$ monoidal.
\end{theorem}

\begin{remark}
\label{rem:solidification_of_solid_is_id}
Using the Yoneda lemma, adjointness, and the full faithfulness of the inclusion $\Mod_{\ul{\mc{A}}} \to \Mod_{\mc{A}}^\Solid$, one sees that the solidification of a solid module is itself, namely, if $M\in \Mod_{\mc{A}}^\Solid$ then $M\otimes_{\ul{\mc{A}}}\mc{A} = M$. A similar remark applies to solid complexes $C \in D(\mc{A})$.
\end{remark}

\begin{example}
\label{ex:solid_modules}
 If $\mc{A}[G]'[S]$ is static for every extremally disconnected $S$ then the category of solid $\mc{A}[G]'$-modules is the same as the category of solid $\mc{A}$-modules with a continuous semilinear $G$-action (\cite[Remark 3.4.4]{Ma22}).
\end{example}

We recall the definition of a nuclear module. 

\begin{definition}[{\cite[Definition 13.10]{CS19an}}]
\label{def:nuclear_modules}
1. Let $M$ be a solid $\mc{A}$-module. We say that $M$ is
nuclear if for all extremally disconnected sets $S$ the natural map
$$\ul{\mt{Hom}}_{\mc{A}}(\mc{M}[S], \mc{A})\otimes_{\mc{A}}M \ra \ul{\mt{Hom}}_{\mc{A}}(\mc{M}[S], M)$$
is an isomorphism.

2. Let $C\in D(\mc{A})$ be a solid $\mc{A}$-complex. We say that $C$ is nuclear if for all extremally disconnected sets $S$ the natural map

$$\ul{\mt{RHom}}_{\mc{A}}(\mc{M}[S], \mc{A})\otimes_{\mc{A}}^LC \ra \ul{\mt{RHom}}_{\mc{A}}(\mc{M}[S], C)$$
is an isomorphism.
\end{definition}

The following is a consequence of \cite[Proposition 5.35]{andreychev2021pseudocoherent}.

\begin{proposition}
Let $M$ be a nuclear $\mc{A}$-module and let $S$ be a profinite set. Then for every complex of solid $\mc{A}$-complex $C\in D(\mc{A})$, the natural map
$$\ul{\RHom}_{\mc{A}}(\mc{M}[S],\mc{A})\otimes^L_{\mc{A}}M \to \ul{\RHom}_{\mc{A}}(\mc{M}[S],M)$$
is an isomorphism.
\end{proposition}

Finally, we recall definitions related to morphisms of analytic rings.

\begin{definition}
\label{def:morphism_of_an_rings}
1. (\cite[$\mathsection{7}$]{CS19con}) A map of analytic
rings $\mc{A}\to \mc{B}$ is the data of a map on the underlying condensed rings which maps $\mc{A}[S]$ to $\mc{B}[S]$ for every extremally disconnected set $S$. Given such a map, there is an induced monoidal base change functor from solid $\mc{A}$ modules to solid $\mc{B}$ modules, denoted $\otimes_{\mc{A}}\mc{B}$. It is given by the composition of $\otimes_{\ul{\mc{A}}}\ul{\mc{B}}$ and $\otimes_{\ul{\mc{B}}}\mc{B}$. Similarly, there is a monoidal functor $\otimes_{\mc{A}}^L\mc{B}$ on the derived category, which is its left derived functor. Given two maps $\mc{A}\to \mc{B}$ and $\mc{A}\to \mc{C}$, the pushout $\mc{B}\otimes^L_{\mc{A}}\mc{C}$ exists, though it might not be static (\cite[Proposition 12.12]{CS19an}).

2. Let $\mc{A}$ be an analytic ring and let $\ul{\mc{B}}$ be a condensed $\ul{\mc{A}}$-algebra. Then \cite[Proposition 12.8]{CS19an} \cite[Definition 2.3.13]{Ma22} describe an induced analytic ring structure induced on $\mc{B}$. 

3. (\cite[Definition 12.13]{CS19an}) A map of analytic rings $\mc{A} \to \mc{B}$ is steady if for all maps $g:\mc{A} \to \mc{C}$, with $\mc{C}$ possibly animated, and for all $M \in D(\mc{C})$, the object $M\otimes^L_{\mc{A}}\mc{B}$, regarded in $D(\mc{B}\otimes^L_{\mc{A}}\mc{C})$, lies in $D(\mc{C})$ when restricted
to $\mc{C}$.

4.  (\cite[Definition 12.16]{CS19an}) A map of analytic rings $\mc{A}\to \mc{B}$ is a localization if the forgetful functor $D_{\ge 0}(\mc{B})\to D_{\ge 0}(\mc{A})$ is fully faithful. It is a steady localization if it is a localization and steady.
\end{definition}

\begin{remark}
\label{rem:induced_an_structures}
    1. When $\mc{A}$ is commutative and $\ul{\mc{B}}$ is a condensed $\ul{\mc{A}}$-algebra, the induced analytic ring structure $\mc{B}$ is defined as the completion of the uncompleted analytic ring whose functor of measures is $S \mapsto \mc{\ul{B}}[S]\otimes_{\mc{\ul{A}}}\mc{A}$. 
    
    2. If $\ul{\mc{B}}$ is a solid $\mc{A}$-module, then the formula of 1 shows that completion is unnecessary, and hence $\mc{B}[S] = \mc{\ul{B}}[S]\otimes_{\mc{\ul{A}}}\mc{A}$.

    3. If $\ul{\mc{B}}$ is a solid $\mc{A}$-module, and so that 2 applies, then a $\ul{\mc{B}}$-module is $\mc{B}$-solid if and only if its restriction to $\ul{\mc{A}}$-modules is $\mc{A}$-solid. To see this, reduce by Remark \ref{rem:solidification_of_solid_is_id} to showing that for a $\ul{\mc{B}}$-module $M$ we have $M\otimes_{\ul{\mc{A}}}\mc{A}=M\otimes_{\ul{\mc{B}}}\mc{B}$. We may then reduce to the case $M=\ul{\mc{B}}[S]$, which follows from 2.
\end{remark}





\section{Solid functional analysis in mixed characteristic}

In this section we introduce the solid spaces which appear in the article. The basic idea is this: we fix a coefficient Banach ring and general Banach or Smith spaces are defined to be orthonormizable modules over it. In the characteristic 0 setting, this coefficient Banach ring can be taken to be $\QQ_p$.

\subsection{Banach pairs}

The following are the same as the $\ZZ_p$-Tate algebras appearing in \cite{Po24}.

\begin{definition}
\label{def:banach_pair}
A Banach pair is a complete Tate Huber pair $(B,B^+)$ together with a morphism $(\ZZ_p, \ZZ_p) \rightarrow (B,B^+)$. 
\end{definition}
We often omit $B^+$ from the notation and simply say that $B$ is a Banach ring when $B^+$ is clear from the context.
 
If $B$ is a Banach ring, we may choose a topologically nilpotent unit $\varpi \in B$. There is\footnote{We emphasize that it could be the case that some rational, non integral power of $\varpi$ lies in $B$, but we still make it so that the valuation is $\ZZ$-valued. Thus if $\varpi^{m/n}\in B^+$ for $m \geq n$ we have $\val_{\varpi}(\varpi^{m/n})=\floor{m/n}.$} a natural $\ZZ$-valued $\varpi$-adic valuation $\val_\varpi$ on $B$ such that $B^+ = B^{\val_\varpi \geq 0}$.  It induces the topology on $B$, and $B^+$ is $\varpi$-adically complete and separated.

\begin{example}
\label{ex:banach_rings}
1. Let $B$ be a Banach $\QQ_p$-algebra with unit ball $B^+$. Then $(B,B^+)$ is a Banach pair with topologically nilpotent unit $p$.

2. Let $B = \FF_p((X))$ and $B^+=\FF_p[[X]]$ taken with their $X$-adic topology. Then $(B,B^+)$ is a Banach pair with topologically nilpotent unit $X$. Similarly, we can take a perfected version of $(B,B^+)$ such as $(\FF_p((X^{1/p^\infty})),\FF_p[[X^{1/p^\infty}]])$.

3. For $1/r\in \ZZ[1/p]_{>0}$, we have the rings  $\widetilde{\mathbf{A}}^{(0,r]} = \A_{\mathrm{inf}}\langle p/[\varpi]^{1/r}\rangle[1/[\varpi]]$ and ${\mathbf{A}}^{(0,r]}_{\QQ_p}=(\ZZ_p[[T]][1/T])_{p}^{\wedge}\cap \widetilde{\mathbf{A}}^{(0,r]}$ for $T=[(1,\zeta_p,\zeta_p^2,...)]-1$, taken with their $[\varpi]$-adic and $T$-adic topology respectively. They are Banach rings\footnote{The valuation of $\widetilde{\mathbf{A}}^{(0,r]}$ appearing in the literature is not the same as ours (for instance, there the valuation of $[\varpi]^{1/r}$ is $1/r$ which is not the case for us). However, these valuations are equivalent.}. These rings appear as coefficient rings in the theory of $(\varphi,\Gamma)$-modules.

4. Given a positive rational number $0 < \lambda = n/m$ with $(n,m) = 1$, set
$$(\mc{O}_{\lambda},\mc{O}_{\lambda}^\circ)=(\ZZ_p[[\varpi]]\angles{p^m/\varpi^n}[1/\varpi],\ZZ_p[[\varpi]]\angles{p^m/\varpi^n}).$$
Here, $\mc{O}_{\lambda}^\circ$ is the $(p,\varpi)$-adic completion of $\ZZ_p[[\varpi]][p^m/\varpi^n]$. It is also equal to its $\varpi$-adic completion since $p^m = (p^m/\varpi^n)\cdot \varpi^n$. We equip it with the $\varpi$-adic valuation $\val_{\varpi}$ such that $\mc{O}_{\lambda}^\circ=(\mc{O}_{\lambda})^{\val_{\varpi}\geq 0}.$ The pair $(\mc{O}_{\lambda},\mc{O}_{\lambda}^\circ)$ is then a Banach pair.
\end{example}

\begin{remark}
\label{rem:universal_banach_ring}
Example 4 of \ref{ex:banach_rings} is in a sense universal, as we shall now explain. If $(A, A^+)$ is a complete Huber pair over $(\ZZ_p, \ZZ_p)$, then we have
$$\Hom_{(\ZZ_p, \ZZ_p)}((\mc{O}_{\lambda}, \mc{O}_{\lambda}^\circ),(A,A^+)) \cong \{f\in A^+, f\text{ is invertible and } |p|\leq |f|^\lambda\}.$$ Thus, a Banach pair is the same as a complete Huber pair over $(\ZZ_p,\ZZ_p)$ with a map from $(\mc{O}_\lambda,\mc{O}_\lambda^\circ)$ for some $\lambda$. We can also form a single space which classifies all Banach pairs: this is the pseudorigid open disc $\bb{D}=\bigcup_{\lambda \in \QQ_{>0}}\bb{D}_{\lambda},$ where $\bb{D}_\lambda=\Spa(\mc{O}_{\lambda},\mc{O}_{\lambda}^\circ)$ is the $\lambda$-elementary pseudorigid closed disc defined by the condition $|p|\leq|\varpi|^\lambda$ (see $\mathsection{4}$ of \cite{Lou17}). By the above, we have
$$\Hom_{\Spa(\ZZ_p, \ZZ_p)}(\Spa(A, A^+),\bb{D})\cong \{f\in A^+, f\text{ is invertible and } |p| < 1\},$$
So that $\bb{D}$ classifies Banach pairs (with a choice of a topologically nilpotent unit).\end{remark}

\begin{definition}
\label{def:slope}
Given a Banach pair, we set
$$\slope(B,B^+) = \sup \{\lambda \in \QQ_{>0}:|p|\leq |\varpi|^\lambda\}.$$
Equivalently, $\slope((B,B^+)) \ge \lambda$ if there exists a map $(\mc{O}_{\lambda},\mc{O}_{\lambda}^\circ)\rightarrow (B,B^+)$ over $(\ZZ_p,\ZZ_p)$.  The set on which the supremum is taken is nonempty by virtue of the previous remark. Geometrically, $\slope(B,B^+)$ is the largest $\lambda$ such that the $\Spa(B,B^+)$-point of $\bb{D}$ lies $\bb{D}_{\lambda}$.
\end{definition}

\begin{remark}
1. Examples 1-3 of \ref{ex:banach_rings} are all of slope $\ge 1$ (in Example 3.2.3, this requires that $r$ is sufficiently small). Throughout this article we will usually restrict to this case. Doing this makes several aspects technically cleaner and this assumption seems to hold in most applications of interest.

2. If $\lambda \le \lambda'$ and $(B,B^+)$ is of slope $\ge \lambda'$ then it is also of slope $\ge \lambda$. To see this, write $\lambda'=m'/n'$ and $\lambda = m/n$, so that $mn'-m'n \ge 0$. We have
$$(p^{m}/\varpi^n)^{n'} = p^{mn'-m'n}\cdot [p^{m'}/\varpi^{n'}]^n \in B^+,$$
and since $B^+$ is integrally closed, we have $p^{m}/\varpi^n \in B^+$. Thus, the supremum is taken over a half open interval starting at 0.
\end{remark}

\subsection{Banach pairs as analytic rings}
\label{subsec:banach_pairs_as_an_rings}

By Example \ref{ex:analytic_rings}.1, one may associate to a Banach pair $(B,B^+)$ an analytic ring $(B, B^+)_\bs$ which has $\ul{B}$ as its underlying condensed ring. For brevity, we denote this analytic ring by $B_{\bs}$, with its functor of measures which to an extremally disconnected set $S$ (and, after performing an extension as in \ref{rem:rem_on_an_rings}, to a profinite set $S$) associates a $\ul{B}$-module $B_\bs[S]$. We then have the associated category of solid $B_{\bs}$-modules. Similarly, we have the complete Huber pair $(B^+,B^+)$ and its associated analytic ring $B^+_\bs$, for which we can consider the functor of measures $S\mapsto B^+_{\bs}[S]$ and solid $B_{\bs}^+$-modules.
Since $\ul{B}=\varinjlim_n \ul{B}^+\cdot \varpi^{-n}$, the analytic ring $B_\bs$ is a solid $B_\bs^+$-module.

It turns out that the following class of Banach pairs is much more nicely behaved:

\begin{definition}
\label{def:residually_of_finite_type}
A Banach pair $(B,B^+)$ is residually of finite type if $B^+/\varpi$ is a finitely generated $\ZZ$-algebra.
\end{definition}

Fortunately, in practice this assumption holds in example for interest. For instance, every Banach pair appearing in Example \ref{ex:banach_rings} is residually of finite type, except for perfected rings such as $\FF_p((X^{1/p^\infty}))$ or $\wt{\bf{A}}^{(0,r]}$, but usually one does not want to take these as the base Banach pair anyway.
For the rest of this section, we assume $(B,B^+)$ is residually of finite type. One significant advantage of such Banach pairs is that their associated functors of measures have a simple formula, as in the next proposition. 

\begin{proposition}
\label{prop:measure_description_of_banach_pairs} Let $S$ be a profinite set and let $I$ be an index set such that \footnote{Such a set exists by \cite[Corollary 5.5]{CS19con}.} $\ZZ_{\bs}[S] \cong \prod_I \ZZ$. Then $B^+_\bs[S]=\prod_I B^+$ and $B_{\bs}[S] = (\prod_IB^+)[1/\varpi]$.
\end{proposition}

\begin{proof}
This comes down to unpacking \cite[Theorem 3.27]{andreychev2021pseudocoherent}, which constructs $(A,A^+)_\bs$ from a complete Huber pair $(A,A^+)$. 

To do this, recall the notion of a quasi-finitely generated module (\cite[$\mathsection{3.1}$]{andreychev2021pseudocoherent}). Given a complete Huber pair $(A,A^+)$, let $(A_0,I)$ be a pair of definition and let $R \subset A_0$ be a finitely generated $\ZZ$-algebra. An $R$-submodule $M$ of $A^+$ is said to be quasi-finitely generated if $M=\varprojlim_n M_n$ for finitely generated $R$-submodules $M_n \subset A/I^n$ satisfying that $M_n \to M_{n-1}$ is surjective. It turns out that for a given $M$, this condition is independent of $I$, and that $M$ is a closed submodule of $A$, namely, for any $k \ge 0$ there exists an $l \ge 0$ such that $I^lM\subset I^k$. The collection of pairs $(R,M)$ where $R \subset A^+$ is a finitely generated $\ZZ$-algebra and $M\subset A$ is a quasi-finitely generated $R$-submodule is a directed poset, and it is shown in loc. cit. that if $\ZZ_{\bs}[S] = \prod_I \ZZ$ then $(A,A^+)_{\bs}[S]=\varinjlim_{R, M }\prod_I \ul{M}$, the colimit taken over the mentioned poset.

Since $B^+/\varpi$ is a finitely generated $\ZZ$-algebra, we may consider the finitely generated $\ZZ$-algebra $R\subset B^+$ generated by lifts of the generators of $B^+/\varpi$ and by $\varpi$. As $B^+$ is $\varpi$-torsionfree, one proves by an elementary argument that for every $k, n \ge 0$, each $\varpi^{-k}B^+/\varpi^nB^+$ is a finitely generated $R$-module. Thus restricting to these $M$ which are contained in $\varpi^{-k} B^+$, the poset has a maximal element which is $\varprojlim_n \varpi^{-k}B^+/\varpi^nB^+=\varpi^{-k}B^+$. It follows that $\Bplusbs[S] =\prod_I \ul{B}^+,$ and (using that every quasi-finitely generated $M\subset B$ is bounded) that
$$\Bbs[S] =\varinjlim_{k}\prod_I \varpi^{-k} \ul{B}^+ = (\prod_I\ul{B}^+)[1/\varpi],$$
as required.\end{proof}

\begin{corollary}

1. The analytic ring structure on $\ul{B}$ induced from $B^+_{\bs}$ coincides with $B_{\bs}$.

2. The map $\Bplusbs\to \Bbs$ of analytic rings is a steady localization.    
\end{corollary}

\begin{proof}
1. The induced analytic ring structure on $\ul{B}$ from $B^+_{\bs}$ is given by mapping an extremally disconnected set $S$ to $\ul{B}[S] \otimes_{\ul{B}^+}B^+_{\bs}$ (\cite[Proposition 12.8]{CS19an}), and 
\begin{align*}
\ul{B}[S] \otimes_{\ul{B}^+}B^+_{\bs} &= \varinjlim_k (\varpi^{-k}\ul{B}^+[S])\otimes_{\ul{B}^+}\Bplusbs\\
&=\varinjlim_k (\varpi^{-k}\ul{B}^+[S]\otimes_{\ul{B}^+}\Bplusbs)\\
&=\varinjlim_k(\varpi^{-k}\cdot \Bplusbs[S])=\Bplusbs[S][1/\varpi],
\end{align*}
which is equal to $\Bbs[S]$ by the previous proposition.

2. We start by showing the map is steady. By part 1 and \cite[Proposition 13.14]{CS19an}, it suffices to show that $\ul{B}$ is a nuclear $\Bplusbs$-module. Indeed, writing $\ul{B}=\varinjlim_k \varpi^{-k}B^+$, each $\varpi^{-k}B^+$ is a compact $\Bplusbs$-module and the inclusion maps $\varpi^{-k}B^+\to \varpi^{-(k+1)}B^+$ are of trace class, so $\ul{B}$ is even basic nuclear.

To show $\Bplusbs\to \Bbs$ is a localization, it suffices by \cite[Exercise 12.17]{CS19an} to show that the natural map of analytic rings $\Bbs \to \Bbs \otimes_{\Bplusbs}\Bbs$ is an isomorphism. Here $\Bbs \otimes_{\Bplusbs}\Bbs$ denotes the pushout of $\Bplusbs \rightarrow B_\bs\leftarrow\Bplusbs$. Since the map $\Bplusbs \rightarrow B_\bs$ is steady, \cite[Proposition 12.14]{CS19an} implies that $\Bbs \otimes_{\Bplusbs}\Bbs$ is the same as the base change of the solid $\Bplusbs$-module $\Bbs$ from $\Bplusbs$ to $\Bbs$, from which the claimed isomorphism easily follows.
\end{proof}

\begin{lemma}

1. Let $S$ be a profinite set. Then there exists an isomorphism of $\Bplusbs$-modules $$\Cont(S,B^+) \cong \widehat{\oplus}_{I} B^+ := \varprojlim_n\oplus_I B^+/\varpi^n$$ for some index set $I$. 

2. Let $I$ be an index set. Then there exists a profinite set $S$ and a retract $\Cont(S,B^+)\to \widehat{\oplus}_I B^+$.
\end{lemma}

\begin{proof}
If $S$ is a profinite set, \cite[Theorem 5.4]{CS19con} shows that the group $\Cont(S,\ZZ)$ is a free abelian group $\oplus_I \ZZ$ for some index set $I$. 

1. Since $\widehat{\oplus}_{I} B^+$ is a $\Bplusbs$-module, it suffices to argue on the level of $\ul{B}^+$-modules. If $A$ is a discrete ring then $\mt{Cont}(S,\bb{Z})\otimes_\bb{Z}A=\mt{Cont}(S,A)$ because $S$ is compact. Hence,
\begin{align*}
\ul{\Cont}(S,B^+) &= \varprojlim_n \ul{\Cont}(S,B^+/\varpi^n)\\
&=\varprojlim_n\oplus_I \ul{B}^+/\varpi^n=\widehat{\oplus}_{I} \ul{B}^+,
\end{align*}
as claimed.

2. Given $I$, we may take $S$ large enough so that $\Cont(S,\ZZ) \cong \oplus_J\ZZ$ with $|I|\le |J|$. By the first part, $\Cont(S,B^+)=\widehat{\oplus}_{J} B^+$, so it retracts to $\widehat{\oplus}_{I} B^+$.
\end{proof}

Given a solid $\Bplusbs$-module $M$ we write $M^\vee:=\RHom_{\Bplusbs}(M,B^+).$

\begin{proposition}
For any index set $I$, we have

1. $(\widehat{\oplus}_IB^+)^\vee = \prod_IB^+$.

2. $(\prod_I B^+)^\vee = \widehat{\oplus}_IB^+$.
\end{proposition}

\begin{proof}
1. We compute:
\begin{align*}
\ul{\RHom}_{\ul{B}^+}(\widehat{\oplus}_I\ul{B}^+,\ul{B}^+) &= \varprojlim_n\ul{\RHom}_{\ul{B}^+/\varpi^n}({\oplus}_I\ul{B}^+/\varpi^n,\ul{B}^+/\varpi^n)\\
&= \varprojlim_n\prod_I\ul{B}^+/\varpi^n
=\prod_I\ul{B}^+.
\end{align*}
2. 
 By \ref{prop:measure_description_of_banach_pairs} we know that $\Bplusbs[S]\cong \prod_JB^+$, with $|J|$ arbitrarily large. It follows that $\prod_I B^+$ is a retract of some $\Bplusbs[S]$, by forgetting the coordinates of $J$ not in $I$. Since $\Bplusbs[S]$ is compact projective, a simple diagram chase shows that $\prod_I B^+$ is projective. Dualizing, and using that $$\RHom_{\ul{B}^+}(\Bplusbs[S],\ul{B}^+)=\RHom_{\ul{B}^+}(\ul{B}^+[S],\ul{B}^+)=\Cont(S,\ul{B}^+)=\widehat{\oplus}_J\ul{B}^+,$$ the retraction also dualizes. Thus $\RHom_{\ul{B}^+}(\prod_I \ul{B}^+,\ul{B}^+)= \Hom_{\ul{B}^+}(\prod_I \ul{B}^+,\ul{B}^+)$, and it is the retraction of $\widehat{\oplus}_J\ul{B}^+$ obtained by forgetting the coordinates of $J$ not in $I$. In other words, it is equal to $\widehat{\oplus}_I\ul{B}^+$. 
\end{proof}

\begin{corollary}
\label{cor:tens_prod_of_prod}
For any two index sets $I,J$ we have $$\prod_I B^+ \otimes^L_{\Bplusbs} \prod_JB^+ = \prod_{I \times J} B^+.$$
\end{corollary}
\begin{proof}
We may adapt the argument of \cite[Proposition 6.3]{CS19con} to our setting. Namely, write $\prod_I B^+$, respectively $\prod_J B^+$, as a retract of $\Bplusbs[S]$, respectively $\Bplusbs[T]$. Then it suffices to show that $\Bplusbs[S]\otimes_{\Bplusbs}^L\Bplusbs[T]=\Bplusbs[S\times T]$, and this holds by \cite[Proposition 2.11]{andreychev2021pseudocoherent}.
\end{proof}

\subsection{Banach and Smith spaces}

In this subsection we fix a Banach pair $(B,B^+)$ which is residually of finite type. We also choose a topologically nilpotent unit $\varpi \in B^+$. It will be evident from what follows that the definitions and results of this subsection are independent of this choice.
\begin{definition}
\label{def:banach_spaces}
1. (i) A classical $B$-Banach space is a topological $B$-module $V$ such that there exists an isomorphism of topological $B$-modules $$V \cong \wh{\oplus}_{I}B := (\varprojlim_s\oplus_{I} B^+/\varpi^s)[1/\varpi]$$ where the right hand side is given its natural $\varpi$-adic topology.

(ii) A classical $B$-Smith space is a topological $B$-module $M$ such that there exists an isomorphism of topological $B$-modules 
$$M \cong (\prod_{I}B^+)[1/\varpi] $$ where $\prod_{I}B^+$ is given the product topology induced from the $\varpi$-adic topology on $B^+$.



2. (i) A solid $B$-Banach space is a $B_\bs$-module $V$ such that there exists an isomorphism of $B_\bs$-modules $$V \cong \wh{\oplus}_{I}B := (\wh{\oplus}_{i\in\in I}B^+)\otimes_{\Bplusbs}\Bbs.$$

(ii) A solid $B$-Smith space is a $\Bbs$-module $M$ such that there exists an isomorphism of $B_\bs$-modules 
$$M \cong (\prod_{I}B^+)\otimes_{\Bplusbs}\Bbs.$$


\end{definition}

For classical and solid $B$-Banach or $B$-Smith spaces, we can define a unit ball, which is the object we get without inverting $\varpi$. For example, for a classical $B$-Banach space it is given by $V^+ = \wh{\oplus}_IB^+$. Of course, this does not only depend on $V$ but also on a specific isomorphism describing $V$. We then have $V=V^+[1/\varpi]$. If $V$ is a classical $B$-Banach space, we furthermore endow it with the $\ZZ$-valued $\varpi$-adic valuation which makes it so that $V^+=V^{\val_\varpi\geq 0}$.

\begin{remark}
One could also think to give the following alternative definition:
a $B$-Banach space (respectively $B$-Smith space) is a (topological/solid) $B$-module $V$ (respectively $B$-module $M$) with a $\varpi$-adically complete (respectively quasicompact) lattice $V^+ \subset V$ with $V^+/\varpi$ discrete (resp $M^+ \subset M$ with $M^+$ quasiseparated). When $(B,B^+) = (\QQ_p,\ZZ_p)$, this alternative definition is equivalent to ours (see \cite[$\mathsection{3}$]{RJRC22}), and this is the approach taken in loc. cit. However, in general our definition (\ref{def:banach_spaces}) is more restrictive. Ultimately, this is because in general $B^+/\varpi$ may not a field and so finding an orthonormal basis for a unit ball satisfying the assumptions of the alternative definition is not always possible. It might be better to regard the definition of the present paper as a preliminary defining only orthonormizable Banach and Smith spaces, which suffices for the purpose of this article. 
\end{remark}

\begin{example}
\label{ex:banach_space}
Let $S$ be a profinite set.  The module of continuous functions from $S$ to $B$, denoted $\mt{Cont}(S,B)$, is a classical $B$-Banach space. Its dual $\mc{M}(S,B)$, the module of $B$-valued measures on $S$, is a classical $B$-Smith space. 
\end{example}
We then have the following theorem, generalizing \cite[Proposition 3.5]{RJRC22}, \cite[Lemma 3.10]{RJRC22} and \cite[Proposition 3.17]{RJRC22}. To prove it, one checks that they can be proven in exactly the same way, as the proof is completely formal given what we have shown.

\begin{theorem}
\label{thm:banach_main_theorem}
1. The functors $W \mapsto \ul{W}, V \mapsto V(*)_{\mt{top}}$ give an equivalence between the category of classical $B$-Banach (respectively classical $B$-Smith) spaces
 and the category of solid $B$-Banach (respectively solid $B$-Smith) spaces. 
 Furthermore, they preserve exact sequences, and the projective tensor product on classical $B$-Banach spaces corresponds to the solid tensor product on solid $B$-Banach spaces.

2. The functor $V \mapsto V^\vee := \ul{\Hom}_B(V,B)$ gives an antiequivalence between the categories of $B$-Banach spaces and $B$-Smith spaces.\footnote{From now on, by part 1 of the theorem, we may omit the adjectives "classical" and "solid" in such statements.} Explictly, we have

(i) $\ul{\mt{Hom}}_{B_\bs^\circ}(\wh{\oplus}_{i \in I}B^+,B^+)= \prod_{i \in I}B^+$ and $\ul{\mt{Hom}}_{B_\bs}(\wh{\oplus}_{i \in I}B,B)= (\prod_{i \in I}B^+)[1/\varpi].$

(ii) $\ul{\mt{Hom}}_{B_\bs^\circ}(\prod_{i \in I}B^+,B^+)= \wh{\oplus}_{i \in I}B^+$ and $\ul{\mt{Hom}}_{B_\bs}((\prod_{i \in I}B^+)[1/\varpi],B)= \wh{\oplus}_{i \in I}B.$


3. Let $V$ be a $B$-Banach space and $W$ a $B$-Smith space. Then
$\ul{\Hom}_{B_\bs}(W,V) = W^\vee\otimes_{B_\bs} V$ and $\ul{\Hom}_{B_\bs}(V,W) = V^\vee\otimes_{B_\bs} W$.
In particular, if $V$ and $V'$ are $B$-Banach spaces (respectively $B$-Smith spaces) then
$(V\otimes_{B_\bs}V')^\vee = V^\vee \otimes_{B_\bs} V'^\vee.$
\end{theorem}

The following two lemma will be useful in $\mathsection{5}$.

\begin{lemma}
\label{lem:banach_spaces_are_nuclear}
Any $B$-Banach space is nuclear as a $B_\bs$-module.
\end{lemma}
\begin{proof}
Repeat the argument appearing in \cite[Corollary 3.16]{RJRC22}.
\end{proof}

\section{Analytic functions and distributions}

Throughout this section let $(B,B^+)$ be a fixed Banach pair, residually of finite type and of slope $\ge 1$. The goal of this section is to give the definitions of analytic functions and distributions of a compact $p$-adic Lie group with coefficients in $B$. These specialize to the familiar objects of characteristic 0 when $(B,B^+) = (\QQ_p, \ZZ_p)$. 

\subsection{Binomial rings and the Amice theorem}
\label{subsection:Binomial rings and the Amice theorem}

In this subsection we introduce rings of analytic functions  which have coefficients in $(B,B^+)$.

We start by understanding the case of $B^+ = \ZZ_p$. It will serve as motivation for the general case. Let $T$ be a variable. It is well known and elementary to prove that the ring $$\mt{Int}=\{f\in \QQ[T]:f(\ZZ)\subset\ZZ\}$$ is freely generated as a $\ZZ$-module by the binomial functions $\{\binom{T}{n}\}_{n\geq 0}$. In particular, for every $n, m \geq 0$ we can write 
$$\binom{T}{n}\cdot\binom{T}{m}=\sum_{k=0}^{n+m}a_k
\cdot \binom{T}{k}$$
where the $a_k$ belong to $\ZZ$. In any ring of binomial functions we introduce from now on, it will always be implicit that we multiply binomial functions by this formula.
Given $d \geq 0$, $\ul{n}\in \ZZ_{\ge 0}^d$ and variables $T_1,...,T_d$ we let $\binom{\ul{T}}{\ul{n}} = \prod_i\binom{T_i}{n_i}$. Given $h \geq 0$, we introduce the following ring of binomial coefficients with a convergence condition: 
$$\ZZ_p^{h\hBin}(\ul{T})=\{ \sum_{\ul{n} \in \ZZ_{\ge 0}^d} b_{\ul{n}} \binom{\ul{T}}{\ul{n}}:b_{\ul{n}} \in \ZZ_p \text{ and } 0 \leq \val_p(b_{\ul{n}}) - \sum_i\val_p(\floor{n_i/p^h}!) \ra \infty \}.$$
 By the Amice theorem (\cite{amice1964interpolation}), there exists an isomorphism between $\ZZ_p^{h\hBin}(\ul{T})$ and the integral functions which are locally analytic on polydiscs of radius $p^{-h}$. This can be restated as follows. For each $\ul{i} \in (\ZZ/p^h\ZZ)^d$ we let $\tilde{\ul{i}}\in \ZZ_p^d$ be an arbitrary lift of $\ul{i}$. Then there exists an isomorphism 
 \begin{align}
\label{amice_isomorphism}
     \ZZ_p^{h\hBin}(\ul{T}) \cong \prod_{\ul{i}\in(\ZZ/p^h\ZZ)^d}\ZZ_p\angles{(\ul{T}+\tilde{\ul{i}})/p^h}.
 \end{align}
We can now define binomial rings with coefficients in a general Banach ring. 
The right hand side of (\ref{amice_isomorphism}) does not make sense in general if $\ZZ_p$ is replaced with $B^+
$, because we are not allowed to divide by $p^h$. However, the left hand side is always sensible, and this motivates us in the following definition. 
\begin{definition} For $h\geq 0$, we define the following functions from $\ZZ_{\ge 0}^d$ to $\ZZ_{\ge 0}$: 
\begin{align*}
v^{h}(\ul{n}) &=\sum_{i=1}^d\val_{p}(\floor{n_i/p^h}!),\\
v_{h}(\ul{n}) &= \floor{|\ul{n}|/p^h(p-1)}.
\end{align*}
\end{definition}

\begin{definition} 
1. We let
\begin{align*}
B^{\Bin,+}(\ul{T})&=\{\sum_{\ul{n} \in \ZZ_{\ge 0}^d} b_{\ul{n}} \cdot \binom{\ul{T}}{\ul{n}}:b_{\ul{n}} \in B \text{, } 0\leq\val_\varpi(b_{\ul{n}})\ra\infty \},\\
B^{h\hBin,+}(\ul{T})&= \{\sum_{\ul{n} \in \ZZ_{\ge 0}^d} b_{\ul{n}} \cdot \binom{\ul{T}}{\ul{n}}:b_{\ul{n}} \in B \text{, } 0\leq\val_\varpi(b_{\ul{n}}) -  v^{h}(\ul{n})\rightarrow\infty \},\\
B_{h\hBin}^+(\ul{T})&= \{\sum_{\ul{n} \in \ZZ_{\ge 0}^d} b_{\ul{n}} \cdot \binom{\ul{T}}{\ul{n}}:b_{\ul{n}} \in B \text{, } 0\leq\val_\varpi(b_{\ul{n}}) -  v_{h}(\ul{n})\rightarrow\infty \}.
\end{align*}
2. We let $B^{\Bin}(\ul{T}):=B^{\Bin,+}(\ul{T})[1/\varpi],$ $B^{h\hBin}(\ul{T}):=B^{h\hBin,+}(\ul{T})[1/\varpi]$ and $B_{h\hBin}(\ul{T}):=B_{h\hBin}^+(\ul{T})[1/\varpi].$
\end{definition}

\begin{example}
Let $(B,B^+) = (\QQ_p,\ZZ_p)$. Then $\QQ_p^{h\hBin,+}(\ul{T}) = \ZZ_p^{h\hBin}(\ul{T})$ and $\QQ_p^{h\hBin}(\ul{T})$ is identified, by the Amice theorem, with the ring of functions on $\ZZ_p^d$  with $\QQ_p$-coefficients which are analytic on discs of radius $p^{-h}$. As $\QQ_p^{h\hBin}(\ul{T})$ is also a $\QQ_p$-Banach space, we see that $(\QQ_p^{h\hBin}(\ul{T}),\ZZ_p^{h\hBin}(\ul{T}))$ is a $\QQ_p$-Banach pair.
\end{example}

\begin{lemma}
\label{lem:binomial_rings_are_banach_spaces}
1. The pairs $(B^{\Bin},B^{\Bin,+})$, 
$(B^{h\hBin},B^{h\hBin,+})$ and $(B_{h\hBin},B_{h\hBin}^+)$ are $B$-Banach spaces.

2. For $h < h'$ we have $B_{h\hBin} \subset B^{h\hBin} \subset B_{h'\hBin}.$
\end{lemma}

\begin{proof}
1. This is explained by the equalities
\begin{align*}B^{\Bin,+}(\ul{T}) &= \wh{\oplus}_{\un{n}\in \ZZ_{\ge 0}^d}B^+\binom{\ul{T}}{\ul{n}},\\
B^{h\hBin,+}(\ul{T}) &= \wh{\oplus}_{\un{n}\in \ZZ_{\ge 0}^d}B^+ \varpi^{v^{h}(\ul{n})}\binom{\ul{T}}{\ul{n}},\\
B_{h\hBin}^+(\ul{T}) &= \wh{\oplus}_{\un{n}\in \ZZ_{\ge 0}^d}B^+ \varpi^{v_{h}(\ul{n})}\binom{\ul{T}}{\ul{n}}.
\end{align*}


2. It is enough to show there exist constants  $c,c' \in \RR
$ depending on $h,h',d$ such that for all $\ul{n}\in \ZZ_{\geq 0}$ we have
$$v_{h'}(\ul{n}) +c' \leq v^{h}(\ul{n}) \leq v_{h}(\ul{n}) + c.$$ To show this, note two things: first, have the classical formula $\val_p(n!)=\frac{n-s_p(n)}{p-1}$, where $s_p$ is the sum of digits of $n$ in base $p$. Second, observe that for $x \in \RR_{\geq 0}$ and $t \in \RR_{>0}$ there is an inequality
\begin{align}
\label{ineq:floor_division}
\floor{x/t} \leq \floor{x}/t < \floor {x/t}+1. 
\end{align}
For the upper bound, we have
\begin{align*}
v^{h}(\ul{n}) =&\frac{1}{p-1}\sum_i(\floor{n_i/p^h}-s_p(n_i/p^h))\\
\leq &\frac{1}{p-1}\sum_i \floor{n_i/p^h}\ \leq \frac{1}{p-1} \floor{\sum_in_i/p^h}\\
\leq &\floor{\sum_in_i/p^h(p-1)}+1 = v_h(\ul{n})+1,
\end{align*}
where in the last $\le$ we have used (\ref{ineq:floor_division}). This shows that we can take $c=1$.
For the lower bound, use $s_p(n) = O(\log(n))$. Using (2) in the first and last inequalities below, we compute that
\begin{align*}
v^{h}(\ul{n}) =&\frac{1}{p-1}\sum_i\floor{n_i/p^h}-O_h(\log(|\ul{n}|))\\
\geq &\sum_i \floor{n_i/(p-1)p^h} -O_h(\log(|\ul{n}|))\\
\geq &\floor{\sum_i n_i/(p-1)p^h} -O_{h,d}(\log(|\ul{n}|))\\
\geq &v_{h'}(\ul{n})\cdot{p^{h'-h}} -O_{h,d}(\log(|\ul{n}|))\\
= &v_{h'}(\ul{n})  + [(p^{h'-h}-1)v_{h'}(\ul{n})-O_{h,d}(\log(|\ul{n}|))],
\end{align*}
and we conclude the proof by observing that since $h'>h$, the term $$(p^{h'-h}-1)v_{h'}(\ul{n})-O_{h,d}(\log(|\ul{n}|))$$ is bounded below by some $c'$ depending on $h,h',d$. 
\end{proof}

The function $v_h$ satisfies a sort of Lipschitz property that will occasionally be useful.

\begin{lemma}
\label{lem:v_h_lipschitz}
If $|\ul{k}| \le |\ul{n}| + |\ul{m}|$ then $v_h(\ul{n})+v_h(\ul{m}) - v_h(\ul{k}) \le |\ul{n}| + |\ul{m}| - |\ul{k}|$.
\end{lemma}
\begin{proof}
It suffices to prove that given $a,b, c \in \ZZ_{\ge 0}$ and $t \ge 1$ with $c\leq a+b$, we have 
\begin{align}
\label{ineq:ineq_to_prove}
a+b-c \geq \floor{a/t} + \floor{b/t} - \floor{c/t}
\end{align}
To prove this, consider the function $f_t(x) = x-\floor{x/t}$. Then one checks that $f_t(a+b) \le f_t(a) +f_t(b)$ (since $\floor{a/t} + \floor{b/t} \le \floor{(a+b)/t}$) and that $f_t(k+1) \ge f_t(k)$ for $k \in \ZZ_{\ge 0}$  (since $t \ge 1$). Hence,
$$f_t(c) \le f_t(a+b) \le f_t(a) + f_t(b),$$
which is equivalent to \ref{ineq:ineq_to_prove}.
\end{proof}

\begin{lemma}
\label{lem:certain_binomial_rings_are_banach_pairs}
The pairs $(B^{\Bin},B^{\Bin,+})$ and $(B^{h\hBin},B^{h\hBin,+})$ are Banach pairs of slope $\ge 1$. 
\end{lemma}
\begin{proof}
This amount to showing $B^{\Bin,+}$ and $B^{h\hBin,+}$ are rings. We may write 
$$\binom{\ul{T}}{\ul{n}}\cdot\binom{\ul{T}}{\ul{m}} = \sum_{\ul{k}\leq\ul{n}+\ul{m}}a_{\ul{k}}\binom{\ul{T}}{\ul{k}}$$ with each $a_{\ul{k}} \in \ZZ$. This immediately implies that $B^{\Bin,+}(\ul{T})$ is a ring. As for $B^{h\hBin,+}(\ul{T})$, suppose that $\val_\varpi(b_{\ul{n}}) \geq v^h(\ul{n})$ and $\val_\varpi(b_{\ul{m}}) \geq v^h(\ul{m})$. We need to show that for $\ul{k} \leq \ul{n} + \ul{m}$ we have
$\val_\varpi(b_{\ul{n}}b_{\ul{m}}a_{\ul{k}}) \geq v^h(\ul{k}),$ for which it suffices to show $\val_{\varpi}(a_{\ul{k}})\geq v^h(\ul{k})-(v^h(\ul{n})+v^h(\ul{m}))$. By the Amice theorem, we know that $$\prod_i\floor{n_i/p^h}!\prod_i\floor{m_i/p^h}!\binom{\ul{T}}{\ul{n}}\binom{\ul{T}}{\ul{m}} \in \ZZ_p^{h\hBin}(\ul{T}),$$ which shows that 
$\val_p(a_{\ul{k}})\geq v^h(\ul{k})-(v^h(\ul{n})+v^h(\ul{m})).$ As $(B,B^+)$ has slope $\ge 1$, we have $\val_\varpi(a_{\ul{k}}) \geq \val_p(a_{\ul{k}}),$ and this allows us to conclude the proof.
\end{proof}

\begin{remark}
The inequality $v_h(\ul{n})+v_h(\ul{m})\geq v_h(\ul{n}+\ul{m})-1$ implies that $B_{h\hBin}$ is a ring. However, $B_{h\hBin}^{+}$ is not a ring in general, even if $(B,B^+)$ has slope $\ge 1$. For example, if $d=1$, $h\in \ZZ_{\ge 0}$, $p\geq3$ and $\varpi=p$ one can check that $\binom{T}{p^h(p-1)}$ belongs to $B_{h\hBin}^{+}$ but its square does not. Consequently, in general $(B_{h\hBin},B_{h\hBin}^{+})$ is not a Banach pair.
\end{remark}

\subsection{Analytic functions and distributions}
\label{subsection:Baker-Campbell-Hausdorff in mixed characteristic}

In this subsection we fix a $d$-dimensional compact $p$-adic Lie group $G$. Let $\mathfrak{g}$ be the Lie algebra of $G$. It is a $\QQ_p$-vector space of dimension $d$ endowed with a Lie bracket operator $[,]$. Recall that sublattices $\mathfrak{g}_0 \subset \mathfrak{g}$ which satisfy $[\mfrk{g}_0,\mfrk{g}_0]\subset p\mfrk{g}_0$ correspond by integration to open uniform subgroups $G_0 \subset G$ (\cite[Theorem 9.10]{DdMS03}). Choose any such sublattice $\mfrk{g}_0$ and choose an identification of $\ZZ_p$-modules $\mfrk{g}_0\cong \ZZ_p^d$. This also gives a homeomorphism $G_0 \cong \ZZ_p^d$ which respects $p$-power subgroups. Let $g_i$ be the element corresponding to the vector which has $1$ at the $i$'th coordinate and $0$ elsewhere. Every element in $G_0$ can be written uniquely as $\ul{g}^{\ul{x}}:=\prod_i g_i^{x_i}$ for an $\ul{x}\in \ZZ_p^d$. 

We may now define analytic functions and distributions.
\begin{definition}
\label{def:spaces_of_funcs_and_dists}
1. We define the following spaces of 
functions:
\begin{align*}
&(\mc{C}(G_0,B),\mc{C}(G_0,B^+)):=(B^\Bin(\ul{T}),B^{\Bin,\circ}(\ul{T})),\\
&(\mc{C}^
{h\han}(G_0,B),\mc{C}^
{h\han}(G_0,B^+)) := (B^{h\hBin}(\ul{T}),B^{h\hBin,\circ}(\ul{T})),\\&(\mc{C}_{h\han}(G_0,B),\mc{C}_{h\han}(G_0,B^+)):=(B_{h\hBin}(\ul{T}),B_{h\hBin}^\circ(\ul{T})),\\
&\mc{C}^{h^+\han}(G_0,B):=\varinjlim_{h<h'} \mc{C}^{h\han}(G_0,B).
\end{align*}

2. We define the following spaces of distributions:
\begin{align*}
&(\mc{D}(G_0,B),\mc{D}(G_0,B^+)):=(\mc{C}(G_0,B)^\vee,\mc{C}(G_0,B^+)^\vee),\\
&(\mc{D}^
{h\han}(G_0,B),\mc{D}^
{h\han}(G_0,B^+)) := (\mc{C}^
{h\han}(G_0,B)^\vee,\mc{C}^
{h\han}(G_0,B^+)^\vee),\\&(\mc{D}_{h\han}(G_0,B),\mc{D}_{h\han}(G_0,B^+)):=(\mc{C}_{h\han}(G_0,B)^\vee,\mc{C}_{h\han}(G_0,B^+)^\vee),\\
&\mc{D}^{h^+\han}(G_0,B) :=\varprojlim_{h<h'} \mc{D}^{h\han}(G_0,B).
\end{align*}
\end{definition}

By Lemma \ref{lem:binomial_rings_are_banach_spaces}, the spaces of functions $\mc{C},\mc{C}^{h \han}$ and $\mc{C}_{h \han}$ are $B$-Banach spaces and the spaces of distributions $\mc{D},\mc{D}^{h \han}$ and $\mc{D}_{h \han}$ are $B$-Smith spaces. By Lemma \ref{lem:certain_binomial_rings_are_banach_pairs}, the pairs of $\mc{C}$ and $\mc{C}^{h \han}$ are Banach pairs of slope $\ge 1$. Note that the limits defining $\mc{C}^{h^+\han}(G_0,B)$ (respectively $\mc{D}^{h^+\han}(G_0,B)$) can be taken over $\mc{C}_{h \han}(G_0,B)$ (respectively $\mc{D}_{h \han}(G_0,B)$) by part 2 of Lemma \ref{lem:binomial_rings_are_banach_spaces}. By Lemma \ref{lem:certain_binomial_rings_are_banach_pairs} and Proposition \ref{prop:ring_structure_on_distribtuions} below, both of $\mc{C}^{h^+\han}(G_0,B)$ and $\mc{D}^{h^+\han}(G_0,B)$ have a natural $B$-algebra structure.

\begin{remark}
\label{rem:bases_of_distributions}
We have $$(\mc{C}(G_0,B),\mc{C}(G_0,B^+)) = (\mt{Cont}(G_0,B), \mt{Cont}(G_0,B^+)).$$ This follows from the mixed characteristic Mahler theorem (\cite[Theorem 2.2]{Po24}). Consequently, $(\mc{D}(G_0,B),\mc{D}(G_0,B^+))$ is identifed with the space of measures; namely,
$$(\mc{D}(G_0,B),\mc{D}(G_0,B^+)) = (B_\bs[G_0],B_{\bs}^+[G_0])$$ is the Iwasawa algebra on $(B,B^+)$. More explictily let $${\bf{c}}^{\ul{n}}:=\prod_i(g_i-1)^{n_i},$$ which is dual to $\binom{\ul{T}}{\ul{n}}$. By part 3 of Theorem \ref{thm:banach_main_theorem}, we have an equality of $\Bplusbs$-modules
$\mc{D}(G_0,B^+) = \prod_{\ul{n}}B^+ {\bf{c}}^{\ul{n}};$ equivalently, 
$$\mc{D}(G_0,B^+) = \{\sum_{\ul{n}\in \ZZ_{\ge 0}^d}b_{\ul{n}}{\bf{c}}^{\ul{n}}:b_{\ul{n}}\in B^+\}.$$ 
Using the same basis, we can use part 4 of Theorem \ref{thm:banach_main_theorem} describe the analytic spaces of distributions as
\begin{align*}
\mc{D}^{h \han}(G_0,B^+) &= \{\sum_{\ul{n} \in \ZZ_{\ge 0}^d} b_{\ul{n}}{\bf{c}}^{\ul{n}}:b_{\ul{n}}\in B, \val_{\varpi}(b_{\ul{n}})\geq -v^h(\ul{n})\},\\
\mc{D}_{h \han}(G_0,B^+) &= \{\sum_{\ul{n} \in \ZZ_{\ge 0}^d} b_{\ul{n}}{\bf{c}}^{\ul{n}}:b_{\ul{n}}\in B, \val_{\varpi}(b_{\ul{n}})\geq -v_h(\ul{n})\}.
\end{align*}
\end{remark}

\begin{proposition}
\label{prop:ring_structure_on_distribtuions}
The pairs $(\mc{D}(G_0,B),\mc{D}(G_0,B^+)$ and $(\mc{D}(G_0,B),\mc{D}_{h \han}(G_0,B^+))$ have a natural structure of (in general noncommutative) $(B,B^+)$-algebras.
\end{proposition}
\begin{proof}
This amounts to giving each of $\mc{D}(G_0,B^+)$ and $\mc{D}_{h \han}(G_0,B^+)$ a $B^+$-algebra structure. Indeed, by the previous remark, $\mc{D}(G_0,B^+)=\Bplusbs[G_0]$ is the space of $B^+$-measures on $G_0$, or what is the same, the Iwasawa algebra of $G_0$ with $B^+$-coefficients. We give $\mc{D}_{h \han}(G_0,B^+)$ the $B^+$-algebra structure by extending that of $\mc{D}(G_0,B^+)$ which is dense inside it. To show this makes sense, write
$${\bf{c}}^{\ul{n}}\cdot{\bf{c}}^{\ul{m}} =\sum_{\ul{k}}b_{\ul{k}}{\bf{c}}^{\ul{k}}$$
in $\ZZ_{p,\bs}[G_0]$.
We need to show that $\val_{\varpi}(b_{\ul{k}}) \geq v_h(\ul{n})+v_h(\ul{m})-v_h(\ul{k})$ to make the formula of multiplication on $\mc{D}_{h \han}(G_0,B^+)$ converge. For $|\ul{k}| > |\ul{n}|+|\ul{m}|$ this is of course automatic. By \cite[Lemma 7.11]{DdMS03}, we have
$${\bf{c}}^{\ul{n}}\cdot{\bf{c}}^{\ul{m}} \in\sum_{|\ul{k}|\leq{|\ul{n}|+|\ul{m}|}} p^{|\ul{n}|+|\ul{m}|-|\ul{k}|}\ZZ_p{\bf}{c}^{\ul{k}} + \sum_{|\ul{k}|>{|\ul{n}|+|\ul{m}|}} \ZZ_p{\bf}{c}^{\ul{k}},$$
which gives $\val_p(b_{\ul{k}}) \geq |\ul{n}|+|\ul{m}|-|\ul{k}|$ for each $\ul{k}$ with $|\ul{k}| \leq |\ul{n}|+|\ul{m}|$. For each such $\ul{k}$, we get
\begin{align*}
  \val_{\varpi}(b_{\ul{k}}) &\geq   \val_p(b_{\ul{k}})\\
  &\geq |\ul{n}|+|\ul{m}|-|\ul{k}|\\
  &\ge v_h(\ul{n})+v_h(\ul{m})-v_h(\ul{k}),
\end{align*}
where the first inequality we used that $(B,B^+)$ is of slope $\ge 1$, and the second inequality follows from Lemma \ref{lem:v_h_lipschitz}.\end{proof}

\begin{remark}
When $(B,B^+) = (\QQ_p,\ZZ_p)$ the product structure on distributions can be defined as the dual of the formal group law $F$ coming from the Baker-Campbell-Hausdorff formula of $\mfrk{g}$, interpreted as a morphism $$F:\mc{C}^{h \han}(G_0,\QQ_p) \rightarrow \mc{C}^{h \han}(G_0,\QQ_p)\times \mc{C}^{h \han}(G_0,\QQ_p),$$
$$\ul{T}\mapsto F(\ul{X},\ul{Y}) =(F_i(\ul{X},\ul{Y}))_{i=1}^d.$$
See \cite[{$\mathsection{5.2}$}]{emerton2017locally} and \cite[Chapter V, {$\mathsection{4}$}]{serre2009lie}. For this morphism to be defined for some $h$, one needs to know the convergence of the formula on a disc of radius $p^{-h}$. This is proved in \cite[Chapter V, {$\mathsection{4}$}, Theorem 2]{serre2009lie}. In our context, we can reverse this logic to obtain a proof of the convergence of the Baker-Campbell-Hausdorff formula. Thus, dualizing the morphism $$\mc{D}_{h \han}(G_0,B)\times \mc{D}_{h \han}(G_0,B) \rightarrow \mc{D}_{h \han}(G_0,B)$$ of Proposition \ref{prop:ring_structure_on_distribtuions} we obtain a morphism $$F:\mc{C}_{h \han}(G_0,\QQ_p) \rightarrow \mc{C}_{h \han}(G_0,\QQ_p) \times \mc{C}_{h \han}(G_0,\QQ_p).$$ When we reinterpert this as a map on binomial rings, we get a map
$$B_{h \hBin}(\ul{T}) \ra B_{h \hBin}(\ul{X},\ul{Y}),$$
which we think of as a mixed characteristic formal group law, or a 
mixed characteristic Baker-Campbell-Hausdorff formula. Its expansion is given in terms of binomial expansions. Namely, for each $\ul{k}$, we have the data of
$$F(\binom{\ul{T}}{\ul{k}})=F_{\ul{k}}(\ul{X},\ul{Y})=\sum_{\ul{n},\ul{m}\in \ZZ_{\ge 0}^d}a_{\ul{n},\ul{m},\ul{k}}\binom{\ul{X}}{\ul{n}}\binom{\ul{Y}}{\ul{m}}\in B_{h \hBin}(\ul{X},\ul{Y}).$$ The association
$$\sum_{\ul{n}}b_{\ul{n}}\binom{\ul{T}}{\ul{n}}\mapsto \sum_{\ul{n}}b_{\ul{n}}F_{\ul{n}}(\ul{X},\ul{Y})$$
then preserves convergence in radius $p^{-h}$. In the case of $(B,B^+) = (\QQ_p,\ZZ_p),$ the $F_{\ul{n}}(\ul{X},\ul{Y})$ are of course determined by the $F_i=F_{\ul{1}_i}(\ul{X},\ul{Y})$ because $\binom{\ul{T}}{\ul{n}}$ are polynomials in $\binom{\ul{T}}{\ul{1}_i}$ over $\QQ_p$. But in general, the $F_{\ul{k}}(\binom{\ul{T}}{\ul{k}})$ are an additional part of the data.
\end{remark}

\begin{example}
Let $G_0$ be the uniform $p$-adic group generated by the two elements $g_1, g_2$ subject to the relation $g_2g_1g_2^{-1} = g_1^p$. We have coordinates on $G_0$ given by matching $\ul{T} = (t_1,t_2) \in \ZZ_p^2$ with $g_1^{t_1}g_2^{t_2}.$ For each $\ul{k}$, the Baker-Campbell-Hausdorff formula gives rise to a power series
$$F_{\ul{k}}(\ul{X},\ul{Y}) = \sum_{\ul{n},\ul{m} \in \ZZ_{\ge 0}^2} a_{\ul{n},\ul{m},\ul{k}}\binom{\ul{X}}{\ul{n}}\binom{\ul{Y}}{\ul{m}},$$
where the coefficient $a_{\ul{n},\ul{m},\ul{k}}$ is exactly the evaluation of the distribution ${\bf{c}}^{\ul{n}}\cdot {\bf{c}}^{\ul{m}}$ on the function $\binom{\ul{T}}{\ul{k}}$; in other words, it is the coefficient of ${\bf{c}}^{\ul{k}}$ in ${\bf{c}}^{\ul{n}}\cdot {\bf{c}}^{\ul{m}}$.
Using the formula
$$(g_2-1)(g_1-1) = -(g_1-1) + \sum_{k=1}^{1+p^2}\binom{1+p^2}{k}(g_1-1)^k+\binom{1+p^2}{k}(g_1-1)^k(g_2-1)$$
we can compute
\begin{align*}
F_{(0,0)}(\ul{X},\ul{Y}) &= 1,\\
F_{(1,0)}(\ul{X},\ul{Y}) &= \binom{\ul{X}}{(1,0)}+\binom{\ul{Y}}{(1,0)}+p^2\binom{\ul{X}}{(0,1)}\binom{\ul{Y}}{(1,0)},\\
F_{(0,1)}(\ul{X},\ul{Y}) &= \binom{\ul{X}}{(0,1)}+\binom{\ul{Y}}{(0,1)},\\
F_{(1,1)}(\ul{X},\ul{Y}) &= \binom{\ul{X}}{(1,0)}\binom{\ul{Y}}{(0,1)} + (1+p^2)\binom{\ul{X}}{(0,1)}\binom{\ul{Y}}{(1,0)} -p^2\binom{\ul{X}}{(0,1)}\binom{\ul{Y}}{(1,1)} ,
\end{align*}
and so on.
Note that this example has a special feature - each of the $F_{\ul{k}}(\ul{X},\ul{Y})$ is a polynomial. This happens in this example only because the commutator $[g_1,g_2]$ of $g_1,g_2$ is polynomial in $g_1,g_2$, which does not have to be true for a general $G_0$.
\end{example}

\subsection{Twisted distribution algebras}
\label{subsec:twisted_distribution_algebras}
Let $G$ and $G_0$ be as in the previous subsection. It is often the case that in applications of interest, the Banach pair $(B,B^+)$ has an action of $G$, and that the $B$-representations under consideration have an action of $G$ which is semilinear over $B$. This happens for example in the theory $(\varphi,\Gamma)$ modules. In this case the representations end up being modules over certain twisted algebras of distributions. 

We now make this precise. As always, we assume our Banach pairs are residually of finite type and of slope $\ge 1$.
\begin{definition}
\label{def:G_Banach_pairs}
1. A $G$-Banach pair is a Banach pair $(B,B^+)$ endowed with a continuous action of $G$ on $B^+$, where $G$ acts by isometries.

2. A $G$-Banach pair $(B,B^+)$ is said to be locally analytic\footnote{More explanation for this naming choice will be given below, see Remark \ref{rem:unambiguity_of_analytic_action}.} for the $G$-action if for some open subgroup $G_0 \subset G$ we have $(g-1)(\varpi^nB^+) \subset \omega^{n+1}B^+$ for all $g \in G_0$ and $n \in \ZZ_{\ge 0}$.
\end{definition}

In particular, we can make any Banach pair into a $G$-Banach pair by letting $G$ acts trivially, and this action is locally analytic. In that case, semilinear representations are the same as linear representations. The definitions and constructions which
follow specialize to those given in the previous section.

\begin{definition}
Let $(B,B^+)$ be a $G$-Banach pair. The $\ul{B}$-algebra $B_{\bs}[G]'$ (respectively $B^+$-algebra $B^+_{\bs}[G]'$) is defined to have the same underlying $\Bbs$-module (respectively $\Bplusbs$-module) structure as same $B_{\bs}[G]$ (respectively $B^+_{\bs}[G]$) and whose product structure is given by 
$g\cdot b=g(b)\cdot g.$ 
\end{definition}

Recall the notation ${\bf{c}}^{\ul{n}}$ of Remark \ref{rem:bases_of_distributions}.
\begin{lemma}
Let $a \in B$ and let $\ul{n}\in \ZZ^d_{\geq 0}$. Then in $B_{\bs}[G]'$, we have $${\bf{c}}^{\ul{n}}\cdot a = \sum_{\ul{k}\leq \ul{n}}a_{\ul{k}}\cdot {\bf{c}}^{\ul{k}}$$ for some $a_{\ul{k}}\in B$ with $\val_{\varpi}(a_{\ul{k}})\geq \val_{\varpi}(a).$
\end{lemma}
\begin{proof}
We prove this by induction on $|\ul{n}|$, the case $|\ul{n}|=0$ being trivial. Let $i$ be minimal such that $n_i > 0$.
As ${\bf{c}}^{\ul{n}} = \prod_{j\geq i}(g_j-1)^{n_j},$ we may use the inductive assumption to write
$${\bf{c}}^{\ul{n}}\cdot a =(g_i-1)\cdot{\bf{c}}^{\ul{n}-1_i}\cdot a = (g_i-1)\cdot(\sum_{\ul{k}\leq \ul{n}-1_i}a_{\ul{k}} \cdot {\bf{c}}^{\ul{k}})$$
with $\val_\varpi(a_{\ul{k}})\geq \val_{\varpi}(a)$ for each $\ul{k} \leq {\ul{n}}-1_i$.
Now writing
$$(g_i-1)\cdot a_{\ul{k}} = g_i(a_{\ul{k}})\cdot(g_i-1)+(g_i(a_{\ul{k}})-a_{\ul{k}})$$
we get $${\bf{c}}^{\ul{n}}\cdot a = \sum_{\ul{k}\leq\ul{n}}(g_i(a_{\ul{k}-1_i})+g_i(a_{\ul{k}})-a_{\ul{k}})\cdot {\bf{c}}^{\ul{k}},$$where we understand that $a_{\ul{k}}=0$ unless $\ul{0}\leq\ul{k}\leq \ul{n}-1_i.$ Observing that $$\val_{\varpi}(g_i(a_{\ul{k}-1_i})+g_i(a_{\ul{k}})-a_{\ul{k}}))\geq\val(a),$$ this concludes the proof.
\end{proof}

\begin{proposition}
There is a unique ring structure on $\mc{D}_{h\han}(G_0,B)$, denoted by  (respectively $\mc{D}_{h\han}(G_0,B^+)$) such that the map of $B_\bs$-modules $$B_{\bs}[G_0]'\ra \mc{D}_{h\han}(G_0,B),$$ 
respectively the map of $\Bplusbs$-modules $$B^{+}_{\bs}[G_0]'\ra \mc{D}_{h\han}(G_0,B^+)$$ is an algebra homomorphism. 
\end{proposition}

\begin{proof}
The map $B_{\bs}[G_0]'\ra \mc{D}_{h\han}(G_0,B)$ has dense image, so if such extensions exist, they are unique. To show the ring structure on $B^+_{\bs}[G_0]'$ extends to $\mc{D}_{h\han}(G_0,B^+)$, it 
suffices to show that if $b, a\in B$ with 
$\val_{\varpi}(b)\geq-v_h({\ul{n}})$ and $\val_{\varpi}(a)\geq-v_h({\ul{m}})$ then in $B^+_{\bs}[G_0]'$ there is an equality
$$(b{\bf{c}}^{\ul{n}})\cdot (a{\bf{c}}^{\ul{m}}) = \sum_{\ul{l}\in\ZZ^d_{\geq 0}}b_{\ul{l}}'{\bf{c}}^{\ul{l}}$$ with $\val_{\varpi}(b_{\ul{l}}') \geq -v_h(\ul{l})$. From Proposition  \ref{prop:ring_structure_on_distribtuions} we know that given $\ul{n},\ul{m}\in \ZZ_{\ge 0}^d$ we have 
$${\bf{c}}^{\ul{n}}\cdot {\bf{c}}^{\ul{m}} = \sum_{\ul{l}} b_{\ul{n},\ul{m},\ul{l}} {\bf{c}}^{\ul{l}}$$ with 
$\val_{\varpi}(b_{\ul{n},\ul{m},\ul{l}})\geq v_h(\ul{n}) + v_h(\ul{m}) - v_h(\ul{l}).$ Thus, using the previous lemma, we have
\begin{align*}
(b{\bf{c}}^{\ul{n}})\cdot (a{\bf{c}}^{\ul{m}}) &=b({\bf{c}}^{\ul{n}}\cdot a){\bf{c}}^{\ul{m}} \\
&=b\cdot(\sum_{\ul{k}\leq\ul{n}}a_{\ul{k}}\cdot {\bf{c}}^{\ul{k}})\cdot {\bf{c}}^{\ul{m}}\\ 
&= \sum_{\ul{k} \leq \ul{n}} ba_{\ul{k}}(\sum_{\ul{l}\in\ZZ_{\ge 0}^d}b_{\ul{k},\ul{m},\ul{l}}{\bf{c}}^{\ul{l}}) \\
&= \sum_{\ul{l}\in\ZZ_{\ge 0}^d}(\sum_{\ul{k}\leq\ul{n}}b\cdot a_{\ul{k}}\cdot b_{\ul{k},\ul{m},\ul{l}} ){\bf{c}}^{\ul{l}},
\end{align*}
with $\val_{\varpi}(a_{\ul{k}})\geq\val_{\varpi}(a)\geq -v_h(\ul{m}).$ We have
\begin{align*}
\val_{\varpi}(b\cdot a_{\ul{k}}\cdot b_{\ul{k},\ul{m},\ul{l}}) &\geq \val_{\varpi}(b)+\val_{\varpi}(a_{\ul{k}})+\val_{\varpi}(b_{\ul{k},\ul{m},\ul{l}})\\
 &\geq-v_h(\ul{l}).
\end{align*}
Hence, setting $b_{\ul{l}}'=\sum_{\ul{k}\leq\ul{n}}b\cdot a_{\ul{k}}\cdot b_{\ul{k},\ul{m},\ul{l}}$ we have $(b{\bf{c}}^{\ul{n}})\cdot (a{\bf{c}}^{\ul{m}}) = \sum_{\ul{l}}b_{\ul{l}}'{\bf{c}}^{\ul{l}}$ with $\val_{\varpi}(b_{\ul{l}}')\geq-v_h(\ul{l}),$ as required.
\end{proof}

We denote the distributions with the twisted algebra structure by $\mc{D}_{h \han}(G_0,B)'$ and $\mc{D}_{h \han}(G_0,B^+)'.$ To each of the condensed rings $$\ul{\mc{A}}=B_{\bs}[G_0]',B^+_{\bs}[G_0]',\mc{D}_{h \han}(G_0,B)',\mc{D}_{h \han}(G_0,B^+)'$$ there is a naturally associated analytic ring $\mc{A}$ whose underlying condensed ring is $\ul{\mc{A}}$ we can give its induced analytic ring structure from either $B_\bs$ or $\Bplusbs$ (see Definition \ref{def:morphism_of_an_rings}.2).

\begin{lemma}
\label{lem:solid_mods_over_B_G_tag}
Solid modules over $B_{\bs}[G_0]'$ (respectively over $\Bplusbs[G_0]'$) are the same as $G$-semilinear solid $B$-representations (respectively solid $B^+$-representations).
\end{lemma}

\begin{proof}
By Example \ref{ex:solid_modules}, it suffices to show that 
$B_{\bs}^+[G_0]'[S]$ and $B_{\bs}[G_0]'[S]$ are static for every extremally disconnected set $S$. We verify this for $B_{\bs}^+[G_0]'[S]$, the verification that $B_{\bs}[G_0]'[S]$ is static is similar. According to Definition \ref{def:morphism_of_an_rings}.2, we have that $B_{\bs}^+[G_0]'[S]$ is equal to the $\Bplusbs$-solidification of $\ul{B_{\bs}^+[G_0]'}[S]$. As a $\ul{B}^+$-module, $\ul{B_{\bs}^+[G_0]'}[S]=\ul{B_{\bs}^+[G_0]'}\otimes_{\ul{B}^+}\ul{B}^+[S]$. Applying $\Bplusbs$-soldification, which is symmetric monoidal, we get an equality of $\Bplusbs$-modules $$B_{\bs}^+[G_0]'[S] =B^+_{\bs}[G_0]'\otimes_{\Bplusbs}\Bplusbs[S],$$ and this is static by Corollary \ref{cor:tens_prod_of_prod}. 
\end{proof}

In particular, $\Dhplus'$ is a solid $\Bplusbs[G_0]'$-module. We wish to clarify the following potential ambiguity in the notation $\otimes^L_{\Bplusbs[G_0]'}\Dhplus'$.
\begin{lemma}
\label{lem:no_ambig_in_base_change}
Let $M \in D(\Bplusbs[G_0]')$. Then the restriction of the base change of $M$ to $\Dhplus'$ to $\Bplusbs[G_0]'$-modules is equal to the $\Bplusbs[G_0]'$-solid tensor product of $M$ with $\Dhplus'$.
\end{lemma}

\begin{proof}
Both of $\Bplusbs[G_0]'$ and $\Dhplus'$ have the induced analytic ring structure from $\Bplusbs$. Hence, according to Remark \ref{rem:induced_an_structures}.3, we have an equality of functors
\begin{align}
\label{eq:tensor_funcs}
\otimes^L_{\ul{\Dhplus'}}\Dhplus' = \otimes^L_{\ul{B}^+}\Bplusbs = \otimes^L_{\ul{\Bplusbs[G_0]'}}\Bplusbs[G_0]'.
\end{align}
The result follows by applying the left and right functors of (\ref{eq:tensor_funcs}) to $$M\otimes^L_{\ul{\Bplusbs[G_0]'}}\ul{\Dhplus'}.$$
\end{proof}

Finally, we have the following lemma which will be used in $\mathsection$\ref{subsec:idempotency_of_dist_algebras}.
\begin{lemma}
\label{lem:tens_prod_of_prod_with_G0}
For any two index sets $I,J$ we have $$\prod_I \Bplusbs[G_0]'\otimes^L_{\Bplusbs[G_0]'} \prod_J\Bplusbs[G_0]' = \prod_{I \times J} \Bplusbs[G_0]'.$$
\end{lemma}
\begin{proof}
Since the base change functor $\otimes_{\Bplusbs}\Bplusbs[G_0]'$ is monoidal, by Corollary \ref{cor:tens_prod_of_prod} it suffices to show that $\prod_IB^+ \otimes_{\Bplusbs} \Bplusbs[G_0]'=\prod_I \Bplusbs[G_0]'.$ When $\prod_IB^+ = \Bplusbs[S]$ for some extremally disconnected $S$, this was shown in the proof of the previous lemma. In the general case, write $\prod_IB^+$ as a retract of some $\Bplusbs[S]$.
\end{proof}





\subsection{Filtrations, gradings and flatness}
\label{subsec:filtrations_grading_flatness}
In this subsection we assume that the $G$-Banach pair $(B,B^+)$ has a locally analytic $G$-action. Let $G_0$ be an open uniform pro-$p$ subgroup such that for $g\in G_0$ and $n \in \ZZ_{\ge 0}$ we have 
\begin{align}
\label{aligned:cong_analyticity_cond}
(g-1)(\varpi^nB^+) \subset \varpi^{n+1}B^+
\end{align}
In particular, note that $G_0$ acts trivially on $B^+/\varpi$. The purpose of this section is to derive some algebraic results on the algebras we have introduced using techniques going back to the seminal paper of Schneider and Teitelbaum (\cite{schneider_teitelbaum2002algebras}).

\begin{lemma}
\label{lem:Bplus_flatness}
The $\Bplusbs$-modules $\Bplusbs[G_0]'$ and $\Dhpluszero'$ are flat over $\Bplusbs$. 
\end{lemma}
\begin{proof}
This is true in general for unit balls in $B$-Smith spaces, because any such unit ball $\prod_I \Bplusbs$ is the $\Bplusbs$-base change of $\prod_I \ZZ_{p,\bs}$, which is flat over $\ZZ_p$.
\end{proof}

Endow $B$ with the $\varpi$-filtration, namely, $\Fil^i(B)=\varpi^iB^+$ for $i \in \ZZ$. Let $I$ be the augmentation ideal in $\Bplusbs[G_0]'$, in other words, the two-sided ideal generated by $\varpi$ and $g-1$ for $g \in G_0$, equivalently the kernel of the ring homomorphism\footnote{Here we have used implicitly that $G_0$ acts trivially on $\Bplusbs/\varpi$, otherwise only the ring homomorphism to the coinvariants of $\Bplusbs/\varpi$ would be defined.} $\Bplusbs[G_0]' \to \Bplusbs/\varpi$ which maps every $g_i-1$ to 0. We endow $\Bplusbs[G_0]'$ with the $I$-adic filtration, i.e. $\Fil^i(\Bplusbs[G_0]')=I^i\Bplusbs[G_0]'$ for $ \ge 0$. We endow $\Bbs[G_0]'=B\otimes_{\Bplusbs}\Bplusbs[G_0]'$ with the tensor product filtration. Similarly, we endow $\Dhplus'$ with the $I_h$-adic filtration, for $I_h$ being the kernel of the homomorphism $\Dhplus'\to \Bplusbs/\varpi$ mapping every $\varpi^{-v_h(\ul{n})}{\bf{c}}^{\ul{n}}$ with $|\ul{n}|\ge 1$ to $0$, and we endow $\Dh' = \Bbs \otimes_{\Bplusbs}\Dhplus'$ with the tensor product filtration. Note that $\Bbs[G_0]'$ and $\Dh'$ are complete for their respective filtrations, since they induce their topologies. The natural map $\Bbs[G_0]'\to \Dh'$ respects the filtrations.

\begin{proposition}
1. The graded ring $\gr(\Bbs[G_0]')$ is isomorphic to a (commutative!) polynomial ring over the ring $(B^+/\varpi)[\pi^{\pm1}]$ ($\pi$ denoting a variable).

2. The map $\Bbs[G_0]' \to \Dhzero'$ induces an isomorphism $\gr(\Bbs[G_0]')\to \gr(\Dhzero').$
\end{proposition}

\begin{proof}
1. Since $G_0$ is uniform, the ring commutators $[g_i-1,g_j-1]$ are 0 mod $p$, hence mod $\varpi\Bplusbs[G_0]'$, since $(B,B^+)$ is of slope $\ge 1$. In addition, (\ref{aligned:cong_analyticity_cond}) implies that each ring commutator $[g_i-1,b]$ is $0$ mod $\varpi\Bplusbs[G_0]'$ for $b \in B^+$. Hence $\gr^\bullet(\Bbs[G_0]') \cong (B^+/\varpi)[\pi^{\pm{1}}][x_1,...,x_d]$ where $x_i$ is the image of $g_i-1$ in $gr^1(\Bbs[G_0]')$.

2. This follows from $\Bbs[G_0]' \to \Dh'$ being a dense injection.
\end{proof}

\begin{corollary}
\label{cor:flatness_of_dists_over_measures}
The rings $\Bbs[G_0]'$ and $\Dh'$ are left and right noetherian and the maps $B\to \Bbs[G_0]'$ and $\Bbs[G_0]' \to \Dh'$ are left and right flat.
\end{corollary}

\begin{proof}
As $(B,B^+)$ is residually of finite type, the ring $B^+/\varpi$ is noetherian, so the the corollary follows from \cite[Propositions 1.1-1.2]{schneider_teitelbaum2002algebras} (which uses results proved in \cite{huishi1996zariskian}).
\end{proof}

\subsection{Functions and distributions on $G$}

The functions and distributions introduced in the previous subsections can be extended to the entire group $G$, rather than just the uniform subgroup $G_0$. Namely, for $D\in \{\mc{D}_{h \han}, \mc{D}^{h \han},\mc{D}^{h^+ \han}\}$, $*\in \{\emptyset,+\}$ and $\bullet\in \{\emptyset,'\}$ we define
$$D(G,B^{*})^{\bullet}=B^{*}_{\bs}[G]^{\bullet}\otimes_{B^{*}_{\bs}[G_0]^{\bullet}}D(G_0,B^{*})^{\bullet}$$
and $C(G,B^{*})^\bullet$ for $C \in \{\mc{C}_{h \han}, \mc{C}^{h \han}, \mc{C}^{h^+ \han}\}$ are defined by dualizing. In addition, we can form 
\begin{align*}
\mc{D}^{\la}(G,B)^{\bullet}&=\varprojlim_h \mc{D}_{h \han}(G,B)^{\bullet}\\ &= \varprojlim_h\mc{D}^{h \han}(G,B)^{\bullet}\\ &= \varprojlim_h\mc{D}^{h^+ \han}(G,B)^{\bullet}
\end{align*}
and
\begin{align*}
    \mc{C}^{\la}(G,B)^{\bullet}&=\varinjlim \mc{C}_{h \han}(G,B)^{\bullet}\\ &= \varinjlim_h\mc{C}^{h \han}(G,B)^{\bullet}\\ &= \varinjlim_h\mc{C}^{h^+ \han}(G,B)^{\bullet}.
\end{align*}
These are both $B$-algebras, where $\mc{D}^{\la}(G,B)^{\bullet}$ inherits its ring structure from the $\mc{D}_{h \han}(G,B)^{\bullet}$ and $\mc{C}^{\la}(G,B)^{\bullet}$ inherits its ring structure from the $\mc{C}^{h \han}(G,B)^{\bullet}.$  As usual, these coincide with the classical constructions when $(B,B^+) = (\QQ_p,\ZZ_p).$

The functorialities are as follows: for $h' > h,$ we have
\begin{align*}
B_{\bs}[G]^\bullet \ra \mc{D}^\la(G,B)^\bullet\ra ...\ra\mc{D}^{h' \han}(G,B)^\bullet\ra \mc{D}_{h \han}(G,B)^\bullet\\\ra \mc{D}^{h^+\han}(G,B)^\bullet
\ra\mc{D}^{h \han}(G,B)^\bullet \ra ...
\end{align*}
and dually
\begin{align*}
\ul{\mt{Cont}}(G,B)^\bullet \leftarrow \mc{C}^\la(G,B)^\bullet\leftarrow ...\leftarrow\mc{C}^{h' \han}(G,B)^\bullet\leftarrow \mc{C}_{h \han}(G,B)^\bullet\\\leftarrow \mc{C}^{h^+\han}(G,B)^\bullet
\leftarrow\mc{C}^{h \han}(G,B)^\bullet \leftarrow ...
\end{align*}

\subsection{Independence of the uniform subgroup}

So far, all of our constructions depend on the choice of the uniform subgroup $G_0$ as well as the choice of basis $g_1,...,g_d \in G_0.$ In this subsection we explain why this choice does not really make much of a difference. Since the results of this subsection will not be used elsewhere in the article, the reader may safely skip forward. Throughout, we let $(B,B^+)$ be as usual a $G$-Banach pair of slope $\ge 1$ with $G$ acting locally analytically. Throughout, we let $I$ denote the augmentation ideal of $\Bplusbs[G_0]'$, i.e. the kernel of the map $\Bplusbs[G_0]' \to \Bplusbs/\varpi$, and we let $J$ be the kernel of $\ZZ_{p,\bs}[G_0]\to \ZZ_p$, so that for each nonzero $\ul{n}$ we have ${\bf{c}}^{\ul{n}} \in J$, but $p \notin J$. Finally, in this section only, write  $u_h(n)=\floor{n/p^h(p-1)}$, so that $v_h(\ul{n}) = u_h(|\ul{n}|)$. Writing ${\bf{c}}^{\ul{k}} =\varpi^{u_h(k)}\varpi^{-u_h(k)}{\bf{c}}^{\ul{k}}$ shows that $$B^+J^k \subset \varpi^{u_h(k)}\Dhpluszero'\subset I^k\Dhpluszero'.$$
 
\begin{lemma}
\label{lem:change_of_uniform_subgroup_easy}
Let $G_0,G_0'\subset G$ be two uniform subgroups with $G_0'\subset G_0$ with respective bases $g_1',...,g_d'$ and $g_1,...,g_d$ determining the spaces of distributions. Then for $h > 0$, there is a canonical inclusion $\mc{D}_{h\han}(G_0',B^+)'\subset \mc{D}_{h\han}(G_0,B^+)'.$
\end{lemma}

\begin{proof}
Let ${\bf{c}}^{\ul{n}}\in \Bplusbs[G_0]',{\bf{c}'}^{\ul{n}}\in \Bplusbs[G_0']'$ be the elements in Remark \ref{rem:bases_of_distributions} which give a basis of $\mc{D}_{h\han}(G_0,B^+)'$ and $\mc{D}_{h\han}(G_0',B^+)'$ respectively. As $G_0'\subset G_0,$ the element ${\bf{c}'}^{\ul{n}}$ lies in the image of the ideal $J^{|\ul{n}|}$ under the map $\ZZ_p[G_0]\ra B^+_{\bs}[G_0]'$. As $(B,B^+)$ is of slope $\ge 1$, one can write ${\bf{c}'}^{\ul{n}} = \sum_{\ul{k}}b_{\ul{k},\ul{n}}{\bf{c}}^{\ul{k}}$ with $\val_{\varpi}(b_{\ul{k},\ul{n}}) \ge \max(0,|\ul{n}|-|\ul{k}|).$ Thus an element $\sum a_{\ul{n}}{\bf{c}'}^{\ul{n}}$ with $\val_{\varpi}(a_{\ul{n}})\geq -v_h(\ul{n})$ can be rewritten as
$$\sum a_{\ul{n}}{\bf{c}'}^{\ul{n}} = \sum_{\ul{k}}(\sum_{\ul{n}}b_{\ul{k},\ul{n}}a_{\ul{n}}){\bf{c}}^{\ul{k}}.$$
Since $$\val_{\varpi}(\sum_{\ul{n}}b_{\ul{k},\ul{n}}a_{\ul{n}})\geq\inf_{\ul{n}}\{\max(0,|\ul{n}|-|\ul{k}|)-v_h(\ul{n})\}\geq -v_h(\ul{k}),$$
(where in the second inequality we have used Lemma \ref{lem:v_h_lipschitz}) this concludes the proof.
\end{proof}

\begin{corollary}
\label{cor:dist_independant_of_basis}
The space of distributions $\mc{D}_{h \han}(G_0,B^+)'$ depends only on $G_0$ and not on the choice of basis $g_1,...,g_d$.
\end{corollary}

\begin{lemma}
\label{lem:commutator_estimates}
Let $G_0$ be a uniform group. 

(i) Let $x,y \in \ZZ_{p,\bs}[G_0].$ The ring commutator $[x,y]$ lies in $(p J, J^p).$

(ii) Let $g_1, g_2 \in G_0$. The ring commutator $[g_1^p, g_2]$ lies in $(p^2J,p J^p, J^{p^2}).$

(iii) Let $g_1, g_2 \in G_0$. The ring commutator $[(g_1-1)^p,g_2-1]$ lies in $(p^2J,J^p, J^{p^2}).$
\end{lemma}

\begin{proof}
(i) It suffices to prove this when $x,y \in G_0$. We have
$$[x,y]=yx(\{x,y\}-1),$$
where $\{x,y\}=x^{-1}y^{-1}xy$ is the group commutator of $x,y$. Since $G_0$ is uniform, $\{x,y\}$ is a $p$-th power in $G_0$. Writing $z^p=\{x,y\}$, and expanding $z^p-1$ in $z-1,$ we get $z^p-1 = \sum_{k=1}^{p}(z-1)^k\binom{p}{k}$ and so $[x,y]\in (pJ,J^p)$. 

(ii) This is similar to (i), except that by \cite[Lemma 2.4]{DdMS03} we know that $[g_1,g_2]$ is a $p^2$-power.

(iii) We have $(g_1-1)^p = (g_1^p-1)+px$ for some $x\in \ZZ_{p,\bs}[G_0].$ We have 
$$[(g_1-1)^p,g_2-1]=[g_1^p-1,g_2-1]+p[x,g_2-1]=[g_1^p,g_2]+p[x,g_2],$$
so we conclude by (i) and (ii). \end{proof}

\begin{lemma}
\label{lem:technical_bound_on_u_h}
Let $h \ge 1$.

1. For $n \ge p$ we have $u_h(n-p) \ge u_h(n)-1.$

2. For $n \ge 1$ we either have $u_h(n) = u_h(n-1)$ or $n \ge 2$, $u_h(n) = u_h(n-1) + 1$ and $u_h(n-1) = u_h(n-2)$.
\end{lemma}

\begin{proof}
1. Using the bound $\floor{x-y} \ge \floor{x} + \floor{-y}$ we have
$$u_h(n-p) \ge u_h(n)+\floor{-1/p^{h-1}(p-1)},$$
and $h \ge 1$ implies $\floor{-1/p^{h-1}(p-1)}  = \floor{-1/(p-1)}=-1.$

2. If $n=1$ then $u_h(n) = \floor{1/p^h(p-1)} = 0 = u_h(n-1)$, so we may assume $n \ge 2$. Now if $u_h(n) \ne u_h(n-1)$ then we must have $u_h(n) = u_h(n-1) + 1$ by Lemma \ref{lem:v_h_lipschitz}. This means that $n$ and $n-1$ lie in consecutive half open intervals of the form $[kp^h(p-1), (k+1)p^h(p-1))$. Since $h \ge 1$, each interval has length $\ge 2$, and so $n-2$ must lie in the same interval as $n-1$, so $u_h(n-2) = u_h(n-1)$.
\end{proof}

For $\un{n}\in \ZZ_{\geq 0}$, let $n_i = pk_i +r_i$ with $0 \leq r_i < p$. We define ${\bf{d}}^{\un{n}}={\bf{c}}^{\un{r}}\cdot {\bf{c}}^{p\un{k}}.$
These elements will only be used in the proof of Corollary \ref{cor:dist_independence_of_uniform_subgroup} below, and can be viewed as auxiliary.


\begin{lemma}
\label{lem:writing_in_d_basis}
Let $G_0$ be a uniform subgroup of $G$ with basis $g_1,...,g_d$. Assume that $h \ge 1$. Then every element of $\mc{D}_{h \han}(G_0,B^+)'$ can be written as a sum $\sum_{\ul{n}}a_{\ul{n}}{\bf{d}}^{\un{n}}$ with $a_{\ul{n}}\in B^+$ and $\val(a_{\ul{n}})\geq -v_h(\ul{n}).$
\end{lemma}

\begin{proof} Recalling the definition of $\mc{D}_{h \han}(G_0,B^+)'$, it is enough to prove the statement of the lemma for the $B^+$-basis elements $\{\varpi^{-u_h(|\ul{n}|)}{\bf{c}}^{\ul{n}}\}_{\ul{n}\in \ZZ_{\ge 0}^d
}$. We can get from ${\bf{c}}^{\ul{n}}$ to ${\bf{d}}^{\ul{n}}$ in $\ZZ_{p,\bs}[G_0]$ by consecutively commuting elements of the form $(g_i-1)^p$ and $g_j-1$ for $i<j.$ Every time we do this, part (iii) of Lemma \ref{lem:commutator_estimates} implies that we produce an element which lies in $(p^2J,p J^p, J^{p^2})\cdot J^{|\ul{n}|-(p+1)}$. Writing $N = |\ul{n}|$ and $J_N:=(p^2J^{N-p},p J^{N-1}, J^{N+p^2-(p+1)})$, we have in $\ZZ_{p, \bs}[G_0]$ that ${\bf{c}}^{\ul{n}} \in {\bf{d}}^{\ul{n}} +J_N.$

Let $$D = \sum_{\ul{n}} B^+\varpi^{-u_h(\ul{n})}{\bf{d}}^{\ul{n}}.$$
 We need to show that each basis element $\varpi^{-u_h(N)}{\bf{c}}^{\ul{n}}$ belongs to $D$, but currently we only know it belongs to $D + \varpi^{-u_h(N)}B^+J_N$. We therefore may reduce to showing that $\varpi^{-u_h(N)}B^+J_N \subset D + I\Dhpluszero$. Indeed, with this given, it follows that  $\varpi^{-u_h(N)}{\bf{c}}^{\ul{n}}$ belongs to 
$D+ I\Dhpluszero$, and since $\mc{D}_{h \han}(G_0,B^+)'$ is $I$-adically complete, this suffices to conclude the proof by successive approximation.

In fact, we shall show that $$B^+J_N =\varpi^{-u_h(N)}B^+(p^2J^{N-p},pJ^{N-1},J^{N+p^2-(p+1)})$$ is contained in $I\Dhpluszero'$. We show this separately for the ideals generating $J_N$.

1. As $h\ge1$ we have by Lemma \ref{lem:technical_bound_on_u_h}.1 that $2+u_h(N-p) \ge 1 + u_h(N)$.
Hence 
\begin{align*}
\varpi^{-u_h(N)}B^+p^2J^{N-p}\subset& \varpi^{-u_h(N)}\cdot \varpi^2\cdot \varpi^{u_h(N-p)}\subset\\
 \varpi\Dhpluszero'\subset &I\Dhpluszero'.
\end{align*}

2. We have 
\begin{align*}
\varpi^{-u_h(N)}B^+J^{N+p^2-(p+1)}\subset& \varpi^{p^2-p-1}\mc{D}_{h \han}(G_0,B^+)'\subset \\ \varpi\Dhpluszero'\subset& I\Dhpluszero'.
\end{align*}
3. We have
$$\varpi^{-u_h(N)}B^+pJ^{N-1} \subset \varpi^{-u_h(N)}B^+\varpi^{1+u_h(N-1)}\Dhpluszero'.$$
Here there are two cases according to Lemma \ref{lem:technical_bound_on_u_h}.2: in the first case, $u_h(N)=u_h(N-1)$ and then $$\varpi^{-u_h(N)}B^+pJ^{N-1}\subset I_B^+\Dhpluszero'.$$ Otherwise, $N\ge 2$ and $u_h(N)=u_h(N-1)+1$ and $u_h(N-1)=u_h(N-2).$ In particular, the first case applies to $B^+pJ^{N-2}$ so that $$B^+pJ^{N-2}\subset  \varpi^{u_h(N-1)+1}\Dhpluszero'$$ and hence $$\varpi^{-u_h(N)}B^+pJ^{N-1}\subset  J\Dhpluszero' \subset I\Dhpluszero,$$ as required.
\end{proof}
\begin{lemma}
\label{lem:commuting_varpi_power_with_c_r}
For each $\ul{n}\in \ZZ_{\geq 0}$ and $h > 0$, the element $\varpi^{-v_h(\ul{n})}{\bf{d}}^{\ul{n}}$ belongs to $\sum_{\ul{0}\leq \ul{t}\leq\ul{r}}{\bf{c}}^{\ul{t}}B^+\varpi^{-v_h(\ul{n})}{\bf{c}}^{p\ul{k}}$.
\end{lemma}

\begin{proof}
Let $t_{\ul{n},h,g} = g^{-1}(\varpi^{-v_h(\ul{n})})/\varpi^{-v_h(\ul{n})}$. It is a unit of $B^+$. Now use
\begin{align*}
\varpi^{-v_h(\ul{n})}(g_i-1)=&(g_i-1)\cdot t_{\ul{n},h,g_i}\cdot \varpi^{-v_h(\ul{n})}  +(t_{\ul{n},h,g_i}-1)\varpi^{-v_h(\ul{n})}\\
\in& (g_i-1)B^+\varpi^{-v_h(\ul{n})} + B^+\varpi^{-v_h(\ul{n})}
\end{align*}
inductively to commute $\varpi^{-v_h(\ul{n})}$ with the $g_i-1$ appearing in ${\bf{c}}^{\ul{r}}$ (this introduces terms of the form ${\bf{c}}^{\ul{t}}$ with $\ul{0} \leq \ul{t}\leq \ul{r}).$
\end{proof}

\begin{proposition}
\label{prop:change_of_uniform_subgroup_hard}
 Let $G_0, G_0'$ be two uniform subgroups with $G_0' \subset G_0.$ Then there exists a constant $c>0$ depending on $G_0, G_0'$ such that for $h$ sufficiently large there is a natural inclusion $$\mc{D}_{h\han}(G_0,B^+)' \subset \Bplusbs[G_0]'\otimes_{\Bplusbs[G_0']'}\mc{D}_{(h-c)\han}(G_0',B^+)'.$$
\end{proposition}

\begin{proof}
For some $t\in \ZZ_{\geq 1}$ large enough we have $G_0^{p^t}\subset G_0'$. By Lemma \ref{lem:change_of_uniform_subgroup_easy}, we may reduce to the case $G_0'=G_0^{p^t}$ and then further to the case $G_0'=G_1:=G_0^p.$ Let $g_1,..,g_d$ be the chosen basis of $G_0$. By Corollary \ref{cor:dist_independant_of_basis}, we may take $g_1^p, ..., g_d^p$ to be the basis of $G_1$.

It now suffices to show that any element of $\mc{D}_{h\han}(G_0,B^+)'$ is congruent to an element of $\Bplusbs[G_0]'\otimes_{\Bplusbs[G_1]'}\mc{D}_{h\han}(G_1,B^+)'$ modulo $\varpi\mc{D}_{h\han}(G_0,B^+)'$. Because of Lemma \ref{lem:writing_in_d_basis} and Lemma \ref{lem:commuting_varpi_power_with_c_r}, we reduce to showing this for elements of the form $\varpi^{-v_h(\ul{n})}{\bf{c}}^{p\ul{k}}$ where $k_i=\floor{n_i/p}$. Write ${\bf{c}}'^{\ul{m}}$ for the usual basis elements in $\Bplusbs[G_1]'$; thus, ${\bf{c}}'^{\ul{m}} = \prod_i(g_i^p-1)^{m_i}$. As $(g_i-1)^p=g_i^p-1$ mod $p$, hence mod $\varpi$, we have $\varpi^{-v_h(\ul{n})}{\bf{c}}^{p\ul{k}}=\varpi^{-v_h(\ul{n})}{\bf{c}'}^{\ul{k}}$ mod $\varpi$. 

By successive approximation, it is now enough to explain why $\varpi^{-v_h(\ul{n})}{\bf{c}'}^{\ul{k}}$ belongs to $\mc{D}_{h-2 \han}(G_1,B^+)'$ for $h$ sufficiently large. To do this, we need to bound the order of the pole of the coefficient $\varpi^{-v_h(\ul{n})}$. The proof will thus be finished provided we show that for $h$ with $p^{h-1}\geq d$ we have $v_h(\ul{n})\leq v_{h-2}(\ul{k}).$ Indeed, we have
\begin{align*}
v_h(\ul{n})&=\floor{|\ul{n}|/p^h(p-1)}\leq\floor{(p|\ul{k}|+d(p-1))/p^h(p-1)}\\&=\floor{(|\ul{k}|+d(p-1)/p)/p^{h-1}(p-1)}.
\end{align*}
Finally, that this is $\le v_{h-2}(\ul{k})=\floor{|\ul{k}|/p^{h-2}(p-1)}$ is guaranteed by the following elementary lemma, applied to $n = |\ul{k}|,t=d(p-1)/p, a = p^{h-1}(p-1)$ and $b = p^{h-2}(p-1)$. \end{proof}

\begin{lemma}
Let $a,b\in \ZZ_{\geq1}$ and let $t \in \RR_{\geq 0}$. Suppose that $b \geq t$ and $a \geq 2b$. Then for all $n \in \ZZ_{\geq 0}$ we have $\floor{(n+t)/a}\leq\floor{n/b}.$
\end{lemma}
\begin{proof}
If $0 \leq n < b$, this is obvious. If $b \leq n < 2b,$ then $$\floor{(n+t)/a}\leq\floor{3b/2b}=1=\floor{n/b}.$$ Finally, if $2b\leq n$ then $$\floor{(n+t)/a}\le(n+t)/a\leq \frac{1}{2}(n/b+1)\leq \frac{1}{2}(\floor{n/b}+2),$$ which is $\le \floor{n/b}$ since $\floor{n/b}\geq 2.$
\end{proof}
\begin{corollary}
\label{cor:dist_independence_of_uniform_subgroup}
The spaces of locally analytic distributions $\mc{D}^\la(G,B^+)'$ and $\mc{D}^\la(G,B)'$ depend only on $G$ and not on the choice of uniform subgroup $G_0$.
\end{corollary}

\begin{proof}
This follows from Lemma \ref{lem:change_of_uniform_subgroup_easy} and Proposition \ref{prop:change_of_uniform_subgroup_hard}. 
\end{proof}

\section{Representations in mixed characteristic}
Let $G$ be a compact $p$-adic Lie group acting on a $G$-Banach pair $(B,B^+)$ which is residually of finite type and of slope $\ge 1$. In this section we define the categories of mixed-characteristic representations of $G$.

\subsection{Continuous representations}
Here we follow the treatment of $\mathsection{4.2}$ of \cite{RJRC22}. We shall consider the category of solid $B_{\bs}[G]'$-modules (respectively solid $B^+_{\bs}[G]'$-modules) as our category of continuous semilinear $B$-representations (respectively $B^+$-representations) of $G$. This is justified by Lemma \ref{lem:solid_mods_over_B_G_tag}.  

As in \cite[{$\mathsection{4.2}$}]{RJRC22} we may extend the definition of continuous functions to complexes.

\begin{definition}
Let $S$ be a profinite set.

1. Let $V\in \mt{Mod}_{B^+_{\bs}}^{\mt{solid}}$. We let
$\ul{\mt{Cont}}(S, V) = \Hom_{B^+_{\bs}}(B^+_{\bs}[S], V ).$

2.  Let $C \in D(B^+_{\bs})$. We let
$\ul{\mt{Cont}}(S, C)= \ul{\mt{RHom}}_{B^+_{\bs}} (B^+_{\bs}[S], C).$\footnote{This is
consistent with part 1 of the definition since $B^+_{\bs}[S]$ is a projective module.}
\end{definition}

If $V\in  \mt{Mod}^{\mt{solid}}_{B^+_{\bs}[G]'}$ (respectively $C \in D(B^+_{\bs}[G]')$) then we can make $\ul{\mt{Cont}}(G, V)$ into an object of $\mt{Mod}^{\mt{solid}}_{B^+_{\bs}[G]'}$ (respectively into an object of $D(B_{\bs}[G]')$) via the conjugation action (denoted $\star_{1,3}$ in \cite[{Proposition 4.25}]{RJRC22}), i.e, via the formula $g(f)(x) = g(f(g^{-1}x)).$ We can make $\RHom_{\Bplusbs[G]'}(\Bplusbs,\ul{\Cont}(G,C))$ into a left $\Bplusbs[G]'$-module via the right regular action of $G$. Namely, if $C$ is a solid $\Bplusbs[G]'$-module, we let $G$ act on  $\eta:\Bplusbs\to \ul{\Cont}(G,C)$ by $(g\eta)(b)(h) = \eta(b)(hg).$

By the same proof of \cite[Proposition 4.25]{RJRC22}, we have 
\begin{lemma}
\label{lem:fixed_elements_after_tensoring_with_funcs}Let $C \in D(B^+_{\bs}[G]')$. There is a natural isomorphism in $D(B^+_{\bs}[G]'):$
$$\ul{\mt{RHom}}_{B^+_{\bs}[G]'}(B^+_{\bs},\ul{\mt{Cont}}(G,C))\xrightarrow{\cong} C.$$
\end{lemma}

\subsection{Analytic representations}

We start by defining analytic vectors. Let us stress that this generalizes the classical constructions in the case of a Banach space over $\QQ_p$ with an action of $G$, as we shall explain below.

\begin{definition}
1. Let $V\in \Mod^{\mt{solid}}_{B_{\bs}[G]'}$. We set
\begin{align*}
V^{h\han} &= \ul{\mt{Hom}}_{B_{\bs}[G]'}(\mc{D}^{h\han}(G,B)',V),\\
V_{h\han} &= \ul{\mt{Hom}}_{B_{\bs}[G]'}(\mc{D}_{h\han}(G,B)',V),\\
V^{h^+\han} &= \varprojlim_{h'>h}V^{h\han}= \varprojlim_{h'>h}V_{h\han}.
\end{align*}
2. Let $C \in D(B_{\bs}[G]').$ We set
\begin{align*}
C^{h\han} &= \ul{\mt{RHom}}_{B_{\bs}[G]'}(\mc{D}^{h\han}(G,B)',V),\\
C_{h\han} &= \ul{\mt{RHom}}_{B_{\bs}[G]'}(\mc{D}_{h\han}(G,B)',V),\\
C^{h^+\han} &= \mt{R}\varprojlim_{h'>h}C^{h\han}=\mt{R}\varprojlim_{h'>h}C_{h\han}.
\end{align*}
In both cases, there is an induced $\Dh'$ or $\mc{D}^{h \han}(G,B)'$ left module structure given by precomposition with multiplication on the right.
\end{definition}


In the case of nuclear modules (Definition \ref{def:nuclear_modules}) we can give a more familiar description of these functors. 
\begin{proposition}
\label{prop:analytic_vectors_for_nuclear}
1. Suppose $M\in\Mod^{\mt{solid}}_{B_{\bs}[G]'}$ is nuclear as a $B_\bs$-module. Then
\begin{align*}
V^{h \han} &= \ul{\Hom}_{B_{\bs}[G]'}(B,\mc{C}^{h \han}(G,B)\otimes_{B_{\bs}}V),\\
V_{h \han} &= \ul{\Hom}_{B_{\bs}[G]'}(B,\mc{C}_{h \han}(G,B)\otimes_{B_{\bs}}V).
\end{align*}
2. Suppose $C\in D(B_{\bs}[G]')$ is nuclear as a solid $B$-complex. Then
\begin{align*}
C^{h \han} &= \ul{\mt{RHom}}_{B_{\bs}[G]'}(B,\mc{C}^{h \han}(G,B)\otimes^L_{B_{\bs}}C),\\
C_{h \han} &= \ul{\mt{RHom}}_{B_{\bs}[G]'}(B,\mc{C}_{h \han}(G,B)\otimes^L_{B_{\bs}}C).
\end{align*}
In both parts, the action of $B_{\bs}[G]'$ on functions from $G$ to $B$ is given by the conjugation action (as in Lemma \ref{lem:fixed_elements_after_tensoring_with_funcs}), the action on the tensor product is the diagonal $G$-action, and the action after applying $\ul{\Hom}_{\Bbs[G]'}$ or $\ul{\RHom}_{\Bbs[G]'}$ is the one induced by right regular action of $G$ (as in the discussion prior to Lemma  \ref{lem:fixed_elements_after_tensoring_with_funcs}).
\end{proposition}

\begin{proof}
We give the proof for 2 and $\mc{C}^{h \han}(G,B)$, proofs in other cases are similar. 
We claim that the natural map
\begin{align}
\label{eq:quasi_iso_from_nuclear}
\mc{C}^{h\han}(G,B)\otimes^L_{B_{\bs}} C\ra\ul{\mt{RHom}}_{B_{\bs}}(\mc{D}^{h \han}(G,B)',C)
\end{align}
is a quasi-isomorphism. Indeed,
since $\mc{D}^h(G,B)'$ is a $B$-Smith space, it is a retract of some solid $B$-module of the form $B_{\bs}[S]$ for $S$ extremally disconnected. Thus the claim reduces to the statement of the map
$$\mt{Cont}(S,B)\otimes^L_{B_{\bs}}C\ra \ul{\mt{RHom}}_{B_{\bs}}(B_{\bs}[S],C)$$
being a quasi-isomorphism, which is exactly the content of $C$ being nuclear.

With the claim given, applying $\ul{\mt{RHom}}_{B_{\bs}[G]'}(B_{\bs},)$ to both sides of the quasi-isomorphism gives the desired isomorphism.
\end{proof}

\begin{remark}

1. By Proposition \ref{lem:banach_spaces_are_nuclear}, our definition of analytic vectors agrees with the usual definitions of analytic vectors via analytic functions when $V$ is a $B$-Banach space.

2. Our definition of analytic vectors generalizes the usual definition of analytic vectors when $(B,B^+)=(\QQ_p,\ZZ_p)$, by \cite[Theorem 4.36]{RJRC22}. In fact, loc. cit. shows the two two possible definitions of analyticity (with functions or with distributions) coincide in this case for general $C$, without any nuclearity assumptions required. We do not know if this is true more generally.
\end{remark}

It follows from the definition and Lemma \ref{lem:fixed_elements_after_tensoring_with_funcs} that given $C \in D(B_{\bs}[G]')$ we have naturals map from $C^{h \han}$, $C_{h \han}$ and $C^{h^+ \han}$ to $C$. We have similar maps for $V \in \Mod^{\mt{solid}}_{B_{\bs}[G]'}.$

The following definition should be compared to \cite[Definition 4.29]{RJRC22}.
\begin{definition}
1. A module $V \in \Mod^{\mt{solid}}_{B_{\bs}[G]'}$ is called upper $h$-analytic (respectively lower $h$-analytic, respectively $h^+$-analytic) if the natural map $V^{h \han}\ra V$ (respectively $V_{h \han}\ra V$, respectively $V^{h^+ \han}\ra V$) is an isomorphism.

2. A complex $C \in D(B_{\bs}[G]')$ is called upper $h$-analytic (respectively lower $h$-analytic, respectively $h^+$-analytic) if the natural map $C^{h \han}\ra C$ (respectively $C_{h \han}\ra C$, respectively $C^{h^+ \han}\ra C$) is a isomorphism.
\end{definition}

\begin{remark}
In the classical setting with $(B,B^+) = (\QQ_p,\ZZ_p)$, an upper $h$-analytic representation is the same as a $G_h$-analytic representation (as in \cite{emerton2017locally}).
\end{remark}

In $\mathsection$\ref{sec: Structural results} below we will provide some structure theorems regarding these categories. It turns out that the notions of lower $h$-analyticity and $h^+$-analyticity are the well behaved ones (as far as we can prove).

\subsection{Locally analytic representations}

In this subsection we define locally analytic vectors and locally analytic representations.

\begin{definition}
1. Let $V \in \Mod^{\mt{solid}}_{B_{\bs}[G]'}$. The space of locally analytic vectors of $V$ is given by
$$V^{\la}=\varinjlim_{h\ra\infty}V^{h \han}=\varinjlim_{h\ra\infty}V_{h \han}=\varinjlim_{h\ra\infty}V^{h^+ \han}.$$
2. Let $C \in D(B_{\bs}[G]')$. The space of locally analytic vectors of $V$ is given by
$$C^{\la}=\varinjlim_{h\ra\infty}C^{h \han}=\varinjlim_{h\ra\infty}C_{h \han}=\varinjlim_{h\ra\infty}C^{h^+ \han}.$$
\end{definition}

\begin{definition}
1. Let $V \in \Mod^{\mt{solid}}_{B_{\bs}[G]'}$. We say that $V$ is locally analytic if the natural map $$V^\la \ra V$$ is an isomorphism.

2. Let $C \in D(B_{\bs}[G]')$. We say that $C$ is locally analytic if the natural map $$C^\la \ra C$$ is a isomorphism.
\end{definition}

\begin{remark}
\label{rem:unambiguity_of_analytic_action}
When $(C,C^+)$ is a $G$-Banach pair over a $G$-Banach pair $(B,B^+)$, there are two possible meanings attached to the statement "$(C,C^+)$ is locally $G$-analytic": either it is meant to be $G$-locally analytic as a Banach pair in the sense of Definition \ref{def:G_Banach_pairs}, or it is $G$-locally analytic as a $B$-Banach space with a $G$-action. We explain why there is no ambiguity and both notions coincide. The point is that for Banach pairs, being $G$-locally analytic is actually intrinsic and does not depend on $(B,B^+)$. This is similar to \cite[Lemma 2.1.5]{camargo2022geometric}. We content with being quite brief and sketchy since this remark will not be used elsewhere in the article. On the one hand, if $(C,C^+)$ is locally analytic as a Banach pair, then for some sufficiently small $G_0 \subset G$ it holds that $(g-1)(\varpi^nC^+)\subset \varpi^{n+1}C^+$ for all $n$. So given $c \in C^+$, and given $g^{\ul{x}}=g_1^{x_1}\cdot...\cdot g_d^{x_d}$, we may write
$$g(c) = \sum_{\ul{n}}\binom{\ul{x}}{\ul{n}}(g_1-1)^{x_1}\cdot...\cdot(g_d-1)^{x_d}(c),$$
so we see that for some sufficiently large $h$, every element of $C$ belongs to $\mc{C}^{h\han}(G_0,C)^{G_0}$, so that $C^{h \han} = C.$ 
Conversely, if $(C,C^+)$ is locally analytic as a Banach space, then arguing as in \cite[Lemma 2.1.5]{camargo2022geometric} the orbit map $C \to \varinjlim_h \mc{C}^{h\han}(G_0,C) $ factors through some $h$ to a map $C \to \mc{C}^{h\han}(G_0,C) = \mc{C}^{h\han}(G_0,B)\otimes_{\Bbs}C$, and even to a map $C^+ \to \mc{C}^{h\han}(G_0,B^+)\otimes_{\Bplusbs} C^+$. Now the action of $G_0^{p^k}$ on $\mc{C}^{h\han}(G_0,B^+)$ maps $\varpi^n\mc{C}^{h\han}(G_0,B^+)$ to $\varpi^{n+1}\mc{C}^{h\han}(G_0,B^+)$, provided we take $k$ large enough (one checks this on the generators of $G_0^{p^k}$, using that $(B,B^+)$ is $G$-locally analytic and of slope $\ge 1$). Hence the same holds for $C^+$.
\end{remark}

\section{Structural results}
\label{sec: Structural results}
As in the previous section, let $(B,B^+)$ be a $G$-Banach pair which is $G$-locally analytic, of slope $\ge 1$, and residually of finite type. Let $G_0$ be an open uniform subgroup such that (\ref{aligned:cong_analyticity_cond}) holds. In this section, we prove structural results about the categories of analytic representations and their cohomologies.

\subsection{The Lazard-Serre and Kohlhaase resolutions}

\begin{lemma} Let $S$ be a profinite set.
Then $B_{\bs}^+ \otimes_{\ZZ_{p,\bs}}^L\ZZ_{p,\bs}[S]= B^+_{\bs}[S].$
\end{lemma}

\begin{proof}
For $S$ extremally disconnected this is true by definition of the base change functor. In general, one can reduce to this case using a retract and Proposition 
\ref{prop:measure_description_of_banach_pairs}.
\end{proof}

Recall the Lazard-Serre resolution (\cite[Théorème 3.2.7, Chapitre V]{lazard1965groupes}). It is a resolution of $\ZZ_{p,\bs}[G_0]$-modules of the trivial module $\ZZ_{p,\bs}$ which has the form
$$0 \ra \ZZ_{p,\bs}[G_0]^{\binom{d}{d}} \ra \ZZ_{p,\bs}[G_0]^{\binom{d}{d-1}} \ra ... \ra\ZZ_{p,\bs}[G_0]^{\binom{d}{0}} \to \ZZ_{p,\bs}\to 0.$$ We denote its differentials by $\partial_{\mt{LS}}^\bullet$. Furthermore, it is equipped with a $\ZZ_{p,\bs}$-linear contracting homotopy $s_{\mt{LS}}^\bullet$ (i.e., a homotopy between the identity and augmentation maps).
Tensor this resolution from the left with $B_{\bs}^+$. The previous lemma implies the following.
\begin{theorem}[Lazard-Serre]
\label{thm:lazard_serre_resolution}
There exists a resolution of $B_{\bs}[G_0]'$-modules
$$0 \ra B^+_{\bs}[G_0]'^{\binom{d}{d}} \ra B^+_{\bs}[G_0]'^{\binom{d}{d-1}} \ra ... \ra B^+_{\bs}[G_0]'^{\binom{d}{0}} \to \Bplusbs \to 0. $$
  It is equipped with a $B_{\bs}^+$-linear contracting homotopy $s_{\mt{LS},B}^\bullet.$
\end{theorem}

\begin{example}
\label{ex:lazard_serre}
If $G_0 = \ZZ_p$ with generator $g$, then the Lazard-Serre complex can be taken to be the resolution
$$0 \ra \Bplusbs[\ZZ_p]' \xrightarrow{g-1} \Bplusbs[\ZZ_p]' \to \Bplusbs \to 0$$
(the $g-1$ multiplication is from the right). We have $s_{\LS,B^+}^{-1}:\Bplusbs \ra \Bplusbs[\ZZ_p]'$ given by the inclusion $\Bplusbs\to \Bplusbs[\ZZ_p]'$ and $s_{\LS,B^+}^{0}:\Bplusbs[\ZZ_p]'\ra \Bplusbs[\ZZ_p]'$. The identity $\partial^{-1} s^0 + s^{-1}\partial^0=\Id$ forces $$s_{\LS,B^+}^0(\sum_{n \ge 0} a_n(g-1)^n)(g-1) = \sum_{n \ge 1} a_n(g-1)^n$$ and thus $$s_{\LS,B^+}^0(\sum_{n \ge 0} a_n(g-1)^n)=\sum_{n \ge 0} a_{n+1}(g-1)^n,$$ in other words, if $f(T) = \sum_{n \ge 0} a_n T^n$, then $$s_{\LS,B^+}^0(f(g-1)) = (f(g-1)-f(0))/(g-1).$$
\end{example}

In the example above, the homotopy divides by the augmentation ideal of $\Bplusbs[G_0]'$ once. This turns out to be a general phenomenon:

\begin{lemma}
\label{lem:lazard_serre_homotopy_division}
Let $I$ be the augmentation ideal of $\Bplusbs[G_0]'$, as in $\mathsection$\ref{subsec:filtrations_grading_flatness}. Then for $n \ge 1$, the homotopy
$s_{\LS,B^+}^\bullet$ maps $I^n B^+_{\bs}[G_0]'^{\binom{d}{\bullet-1}}$ to $I^{n-1}B^+_{\bs}[G_0]'^{\binom{d}{\bullet}}$.
\end{lemma}
\begin{proof}
Using that $s_{\LS,B^+}^\bullet$ is $\Bplusbs$-linear, and that $(B,B^+)$ is of slope $\le 1$, we may reduce to the case $(B,B^+) = (\QQ_p,\ZZ_p)$. 
In this case, Lazard constructs the homotopy $s_{\LS}$ by lifting a homotopy $\bar{s}_{\LS}$ from the Koszul complex of $$\gr^\bullet(\ZZ_{p,\bs}[G_0])=\FF_p[\pi][x_1,...,x_d],$$ where $\gr(s_{\LS})=\bar{s}_{\LS}$ (\cite[Chap. V, 2.1.1]{lazard1965groupes}). We therefore reduce to proving the claim for the Koszul complex $\FF_p[\pi][x_1,...,x_d]$, its homotopy $\bar{s}_{\LS}$ and the ideal $\bar{I}=(\pi,x_1,...,x_d)$. In this case, the Koszul complex is defined inductively from the case $d=1$, and the inductive formulas for $\bar{s}_{\LS}$ (\cite[Chap. V, 1.3.2.2, 1.3.2.3]{lazard1965groupes}) allows one to reduce to the case $d=1$, in which case it is easy to check (take $\gr^\bullet$ of Example \ref{ex:lazard_serre} when $B^+ = \ZZ_p$).
\end{proof}

We will need a version of the Lazard-Serre complex with the $B^+_{\bs}[G_0]'$ replaced by $\mc{D}_{h \han}(G_0,B^+)'$. The following is a generalization of Kohlhaase's resolution (\cite[Theorem 4.4]{kohlhaase2011cohomology}, \cite[Theorem 5.8]{RJRC22}).

\begin{proposition}
\label{prop:integral_kohlhaase_complex}
The Lazard-Serre resolution extends to a complex  of $\mc{D}_{h \han}(G_0,B^+)'$-modules
$$C_{\Koh,B^+} = [0 \ra \mc{D}_{h \han}(G_0,B^+)'^{\binom{d}{d}} \ra \mc{D}_{h\han}(G_0,B^+)'^{\binom{d}{d-1}} \ra $$$$ ...\ra \mc{D}_{h \han}(G_0,B^+)'^{\binom{d}{0}} \to \Bplusbs \to 0].$$
\end{proposition}
\begin{proof}
We are going to extend both the differentials and the homotopy from the Lazard-Serre resolution. We have differentials
$$\partial_{\mt{LS},B^+}^\bullet: B_{\bs}^+[G_0]'^{\binom{d}{\bullet}} \ra B_{\bs}^+[G_0]'^{\binom{d}{\bullet-1}},$$
with $\partial_{\mt{LS},B^+} = 1 \otimes \partial_{\mt{LS}}.$
We are going to explain how it extends continuously to a differential
$$\partial^\bullet_{\mt{Koh},B}:\mc{D}_{h \han}(G_0,B^+)'^{\binom{d}{\bullet}} \ra \mc{D}_{h \han}(G_0,B^+)'^{\binom{d}{\bullet-1}}.$$
By choosing a basis $\{e_i\}$ for $\ZZ_{p,\bs}[G_0]^{\binom{d}{\bullet}}$ and $\{f_j\}$ for $\ZZ_{p,\bs}[G_0]^{\binom{d}{\bullet-1}}$ we reduce to extending continuously the morphism
$$\partial_{\mt{LS},B,ij}^\bullet:B^+_{\bs}\otimes_{\ZZ_{p,\bs}}\ZZ_{p,\bs}[G_0]\cdot e_i\ra B^+_{\bs}\otimes_{\ZZ_{p,\bs}}\ZZ_{p,\bs}[G_0]\cdot f_j.$$
Identifying the domain and codomain with $B^+_{\bs}\otimes_{\ZZ_{p,\bs}}\ZZ_{p,\bs}[G_0],$ this map is given by extending to sums the formula
$$\partial_{\mt{LS},B,ij}^\bullet(b\otimes {\bf{c}}^{\ul{n}}) = b \otimes \partial_{\mt{LS},ij}^\bullet({\bf{c}}^{\ul{n}}) = (b\otimes {\bf{c}}^{\ul{n}})\partial_{\mt{LS},ij}^\bullet(1).$$
In other words, $\partial_{\mt{LS},B,ij}^\bullet$ is given by right multiplication with $\partial_{\mt{LS},ij}^\bullet(1) \in \ZZ_{p,\bs}[G_0]$. Thus, it extends continuously to a left $\mc{D}_{h \han}(G_0,B^+)'$-linear map $$\partial_{\mt{Koh},B^+,ij}^\bullet:\mc{D}_{h \han}(G_0,B^+)'\cdot e_i\ra \mc{D}_{h \han}(G_0,B^+)'\cdot f_j,$$ and so we get the desired maps $\partial_{\mt{Koh},B^+}^\bullet.$ Furthermore, the identities $\partial_{\mt{Koh},B^+}^\bullet\circ\partial_{\mt{Koh},B^+}^{\bullet-1}= 0$ are still satisfied because the inclusion $B_{\bs}[G_0]' \ra \mc{D}_{h \han}(G_0,B)'$ is dense.
\end{proof}

\begin{proposition}
\label{prop:extending_homotopy_to_kohlhaase}
The contracting homotopy $s_{\LS,B^+}^\bullet$ extends to a homotopy $$s_{\Koh,B^+}^{\bullet}:C_{\Koh,B^+}^\bullet\to \varpi^{-1}C_{\Koh,B^+}^{\bullet+1}.$$
\end{proposition}

\begin{proof}
 Again we can reduce to extending a morphism $s_{\mt{LS},B,ij}^\bullet$ on a one-dimensional $B^+_{\bs}[G_0]'$-module to $\mc{D}_{h \han}(G_0,B^+)'$, and again by density of the inclusion $B_{\bs}[G_0]'\ra \mc{D}_{h \han}(G_0,B)'$, there is at most one way to perform this extension. If such an extension exists, it has to be a contracting homotopy, by density. Thus, everything comes down to showing that  that given
$\sum_{\ul{n}}b_{\ul{n}}{\bf{c}}^{\ul{n}}\in \mc{D}_{h \han}(G_0,B^+)',$ the sum
\begin{align}
\label{sum:sum_to_converge}
\sum_{\ul{n}} b_{\ul{n}}s^\bullet_{\mt{LS},ij}({\bf{c}}^{\ul{n}})
\end{align}
converges in $\varpi^{-1}\mc{D}_{h \han}(G_0,B^+)'$. To show this, note that by Lemma \ref{lem:lazard_serre_homotopy_division} we know that $s_{\mt{LS},B}^{\bullet}$ maps things divisible by $I^k$ to things divisible by $I^{k-1}$. This means that  we have $s_{\mt{LS},ij}^\bullet({\bf{c}}^{\ul{n}})=\sum_{\ul{k}}a_{\ul{k}}{\bf{c}}^{\ul{k}}$ such that for $\ul{k}\le \ul{n}$ we have $\val_p(a_{\ul{k}})\geq |\ul{n}|-|\ul{k}|-1$. For $\ul{k}\leq \ul{n}$, we have $|\ul{n}|-|\ul{k}| \ge v_h(\ul{n})-v_h(\ul{k})$, so it follows that
\begin{align*}
\val_{\varpi}(s_{\mt{LS},ij}^\bullet({\bf{c}}^{\ul{n}}))&=\val_{\varpi}(\sum_{\ul{k}}a_{\ul{k}}{\bf{c}}^{\ul{k}})\\ &\geq \inf_{\ul{k}}(v_h(\ul{k})+|\ul{n}|-|\ul{k}|-1)\\&= v_h(\ul{n})-1=\val_{\varpi}({\bf{c}}^{\ul{n}})-1,
\end{align*}$$$$
which establishes the desired convergence of (\ref{sum:sum_to_converge}).
\end{proof}

\begin{corollary}
\label{cor:kolhaase_resolution}
The complex $C_{\Koh,B}:=C_{\Koh,B^+}[1/\varpi]$ has a contracting homotopy. In particular, we get a resolution
$$C_{\Koh,B} = [0 \ra \mc{D}_{h \han}(G_0,B)'^{\binom{d}{d}} \ra \mc{D}_{h\han}(G_0,B)'^{\binom{d}{d-1}} \ra $$$$ ...\ra \mc{D}_{h \han}(G_0,B)'^{\binom{d}{0}} \to \Bbs \to 0].$$
\end{corollary}

\begin{corollary}
\label{cor:dist_tens_triv}
We have\footnote{Recall that by Lemma \ref{lem:no_ambig_in_base_change} there is no ambiguity in the expression $\mc{D}_{h \han}(G_0,B)'\otimes^L_{B_{\bs}[G_0]'}\Bbs$.} $$\mc{D}_{h \han}(G_0,B)'\otimes^L_{B_{\bs}[G_0]'}\Bbs=\Bbs.$$
\end{corollary}
\begin{proof}
We compute:
\begin{align*}
\mc{D}_{h \han}(G_0,B)'\otimes^L_{B_{\bs}[G_0]'}\Bbs &= \mc{D}_{h \han}(G_0,B)'\otimes^L_{B_{\bs}[G_0]'} [...\ra B_{\bs}[G_0]'^{\binom{d}{\bullet}}\ra ...]\\
&= [...\ra \mc{D}_{h \han}(G_0,B)'^{\binom{d}{\bullet}}\ra ...] = \Bbs,
\end{align*}
where the first equality used the Lazard-Serre resolution (Theorem \ref{thm:lazard_serre_resolution}) and the third equality used the Kohlhaase resolution (Corollary \ref{cor:kolhaase_resolution}).
\end{proof}

The following generalization of Lazard's comparison between continuous and analytic cohomology follows immediately.

\begin{theorem}
\label{thm:lazard_comparison}
Let $C\in D(\Bbs[G_0]')$ be lower $h$-analytic. Then 
$$\mt{RHom}_{\Bbs[G]'}(B,C)=\mt{RHom}_{\mc{D}_{h \han}(G,B)'}(B,C).$$
\end{theorem}

\subsection{Idempotency of distribution algebras}
\label{subsec:idempotency_of_dist_algebras}


The goal of this subsection is to prove the following theorem.

\begin{theorem}
\label{thm:idempotent_dists}
We have
\begin{align*}
\mc{D}_{h \han}(G,B)'\otimes^L_{B_{\bs}[G]'}\mc{D}_{h \han}(G,B)' &= \mc{D}_{h \han}(G,B)',\\
\mc{D}^{h^+ \han}(G,B)'\otimes^L_{B_{\bs}[G]'}\mc{D}^{h^+ \han}(G,B)' &= \mc{D}^{h^+ \han}(G,B)'.
\end{align*}
\end{theorem}

\begin{remark}
One can also show (similarly to \cite[Corollary 5.11]{RJRC22}) that 
\begin{align*}
\mc{D}^\la(G,B)\otimes^L_{B_{\bs}[G]}\mc{D}^\la(G,B) &= \mc{D}^\la(G,B).
\end{align*}
We omit the details since this identity will not be used anywhere in the article. We do not know if a similar identity holds for $\mc{D}^\la(G,B)'$ - this seems to be an interesting problem. 
\end{remark}

To prove the theorem, we may reduce to the case $G = G_0$ and treat only the case of $\Dhpluszero$. The proof of it is somewhat tricky, and we proceed in several steps. The idea is to prove the theorem first for the non twisted distributions, and then deduce it for the twisted distributions by using graded techniques.

For this subsection only, we write $\mc{D}_{h,B^+} := \Dhpluszero$ and $\mc{D}_{h,B} = \Dhzero$ to lighten the notation.

Let $C_{\Koh,B^+}$ be the integral Kohlhaase complex (Proposition \ref{prop:integral_kohlhaase_complex}) in the linear setting (so that $\mc{D}_{h,B^+}' = \mc{D}_{h,B^+}$), and let $\tilde{C}_{\Koh,B^+}$ be the same complex but without the last term, so that $$C_{\Koh,B^+}= [... \to \mc{D}_{h,B^+}^{\binom{d}{1}}\to \mc{D}_{h,B^+}^{\binom{d}{0}} \to \Bplusbs \to 0]$$
and 
$$\tilde{C}_{\Koh,B^+}=  [... \to \mc{D}_{h,B^+}^{\binom{d}{1}}\to \mc{D}_{h,B^+}^{\binom{d}{0}} \to 0].$$
\begin{lemma}
\label{lem:d_tens_d_equal_kohl_tens_d}
We have $$\mc{D}_{h,B^+} \otimes_{\Bplusbs[G_0]}^L \mc{D}_{h,B^+} = \tilde{C}_{\Koh,B^+}\otimes_{\Bplusbs}^L\mc{D}_{h,B^+}.$$
\end{lemma}

\begin{proof}
This is essentially shown in the proof of\footnote{Here, one uses strongly that $B^+$ is central in $\Bplusbs[G_0]$. The proof does not seem to work for $\Bplusbs[G_0]'.$} \cite[Proposition 5.10]{RJRC22}. We quickly recall the proof in our context. 

First, we let $\mc{D}_{h,B^+,0}$ denote $\mc{D}_{h,B^+}$ endowed with the trivial action of $G_0$ from the left. If $\iota$ denotes the antipodal map from $\Bplusbs[G_0]$ to itself, then one has an isomorphism of $\Bplusbs[G_0]$-modules
\begin{align}
\label{map:twist_measures_tens_dists}
\phi:\Bplusbs[G_0]\otimes_{\Bplusbs}\mc{D}_{h,B^+}\cong \Bplusbs[G_0]\otimes_{\Bplusbs}\mc{D}_{h,B^+,0}
\end{align}
given as the composition of $1\otimes\iota \otimes 1$ and $1\otimes m$, where $m:\Bplusbs[G_0]\otimes_{\Bplusbs} \mc{D}_{h,B^+}\to \mc{D}_{h,B^+} $ is the multiplication map. In other words, if $g \in G_0
$ and $\mu \in \mc{D}_{h,B^+}$, then $g \otimes \mu$ is mapped to $g \otimes g^{-1}\mu$. To see that $\phi$ is an isomorphism, one checks its inverse maps $g \otimes \mu$ to $g \otimes g\mu$.
Now, $\iota$ extends to $\mc{D}_{h,B^+}$, and so (\ref{map:twist_measures_tens_dists}) extends to 
\begin{align}
\label{map:twist_dists_tens_dists}
\phi:\mc{D}_{h,B^+}\otimes_{\Bplusbs}\mc{D}_{h,B^+}\cong \mc{D}_{h,B^+}\otimes_{\Bplusbs}\mc{D}_{h,B^+,0}.
\end{align}

Let $\wt{C}_{\LS,B^+}$ be the Lazard-Serre complex without the $\Bplusbs$-term, so that $\wt{C}_{\LS,B^+}$ is quasi isomorphic to $\Bplusbs$ and $\mc{D}_{h,B^+}\otimes^L_{\Bplusbs[G_0]}\wt{C}_{\LS,B^+}=\wt{C}_{\Koh,B^+}$.  We then can compute:
\begin{align*}
\mc{D}_{h,B^+}\otimes^L_{\Bplusbs[G_0]}\mc{D}_{h,B^+}&=
 \mc{D}_{h,B^+}\otimes^L_{\Bplusbs[G_0]}(C_{\LS,B^+}\otimes^L_{\Bplusbs}\mc{D}_{h,B^+})\\
 &=\mc{D}_{h,B^+}\otimes^L_{\Bplusbs[G_0]}(C_{\LS,B^+}\otimes^L_{\Bplusbs}\mc{D}_{h,B^+,0})\\
&=(\mc{D}_{h,B^+}\otimes^L_{\Bplusbs[G_0]}C_{\LS,B^+})\otimes^L_{\Bplusbs}\mc{D}_{h,B^+,0}\\
&=\wt{C}_{\Koh,B^+}\otimes^L_{\Bplusbs}\mc{D}_{h,B^+,0}\\
&=\wt{C}_{\Koh,B^+}\otimes^L_{\Bplusbs}\mc{D}_{h,B^+},
\end{align*}

where in the 2nd, respectively 5th equalities we used the isomorphism (\ref{map:twist_measures_tens_dists}), respectively the isomorphism (\ref{map:twist_dists_tens_dists}).
\end{proof}

\begin{lemma}
Let $C = [... \to C_1 \to C_0 \to C_{-1} \to 0]$ be a complex of $\Bplusbs$-modules and let $\tilde{C} = [ ... \to C_1 \to C_0 \to 0]$ be the same complex without $C_{-1}$. Suppose that $H_1(C)$ is killed by $\varpi$ and that $C_0 \to C_{-1}$ is surjective. Then $\ker (H_0(\tilde{C})\to C_{-1})$ is killed by $\varpi$.
\end{lemma}

\begin{proof}
We have a commutative diagram with exact rows
\[
\begin{tikzcd}
  0 \arrow[r]&  \mathrm{ker}(C_0 \to C_{-1}) \arrow[d] \arrow[r]&  C_0 \arrow[r] \arrow[d] & C_{-1} \arrow[d] \arrow[r] &  0 \\
  0 \arrow[r] & K \arrow[r] & H_0(\tilde{C}) \arrow[r] & C_{-1} \arrow[r]&  0
\end{tikzcd}
\]
Let $x \in K$. Then $x \in H_0(\tilde{C}),$ and we can lift it to $\widetilde{x} \in C_0$. By the commutativity of the diagram, it belongs to $\mathrm{ker}(C_0 \to C_{-1})$. As $H_1(C) = \mathrm{ker}(C_0 \to C_{-1})/\mathrm{im}(C_1\to C_0)$ is killed by $\varpi$, we have that $\varpi\tilde{x} \in \mathrm{im}(C_1\to C_0)$; hence, its image $\varpi x$ in $K$ maps to $0$ in $H_0(\tilde{C}) = C_0/im(C_0\to C_1).$
\end{proof}

\begin{proposition}
The map $$H_0(\varpi\mc{D}_{h,B^+}\otimes^{L}_{\Bplusbs[G_0]}\mc{D}_{h,B^+}) \to \varpi \mc{D}_{h,B^+}$$ is an isomorphism.
\end{proposition}

\begin{proof}
Let $C_{\Koh,B^+}$ be the integral Kohlhaase complex (Proposition \ref{prop:integral_kohlhaase_complex}):  
$$C_{\Koh,B^+}= [... \to \mc{D}_{h,B^+}^{\binom{d}{1}}\to \mc{D}_{h,B^+}^{\binom{d}{0}} \to \Bplusbs \to 0]$$
and also let 
$$\tilde{C}_{\Koh,B^+}=  [... \to \mc{D}_{h,B^+}^{\binom{d}{1}}\to \mc{D}_{h,B^+}^{\binom{d}{0}} \to 0].$$
We know by Proposition \ref{prop:extending_homotopy_to_kohlhaase} that the inclusion map $C_{\Koh,B^+} \to \varpi^{-1}C_{\Koh,B^+}$ is homotopic to the zero map. Hence it induces zero on cohomology; so by composing with multiplication by $\varpi$, it follows that $\mt{mul}_{\varpi}:C_{\Koh,B^+}\to C_{\Koh,B^+}$ is the zero map on homologies. Since $\mc{D}_{h,B^+}$ is flat over $\Bplusbs$ (Remark \ref{lem:Bplus_flatness}), $\mt{mul}_{\varpi}$ also induces zero on the homologies of $\tilde{C}_{\Koh,B^+}\otimes_{\Bplusbs}^L\mc{D}_{h,B^+}$.
 Applying the previous lemma for $$\tilde{C}=\tilde{C}_{\Koh,B^+}\otimes^L_{\Bplusbs}\mc{D}_{h,B^+}$$ and $$C=C_{\Koh,B^+}\otimes^L_{\Bplusbs}\mc{D}_{h,B^+},$$ where $C_{-1}  = \Bplusbs \otimes_{\Bplusbs}\mc{D}_{h,B^+} = \mc{D}_{h,B^+},$ we get from Proposition \ref{lem:d_tens_d_equal_kohl_tens_d} that the kernel of
 \begin{align*}
     H_0(\mc{D}_{h,B^+} \otimes^L_{\Bplusbs[G_0]}\mc{D}_{h,B^+})=
     H_0(\tilde{C}_{\Koh,B^+}\otimes^L_{\Bplusbs}\mc{D}_{h,B^+}
) &\to \mc{D}_{h,B^+}
 \end{align*} is killed by $\varpi$. Hence, the map $$H_0(\varpi\mc{D}_{h,B^+}\otimes^{L}_{\Bplusbs[G_0]}\mc{D}_{h,B^+}) \to \varpi \mc{D}_{h,B^+}$$ is an isomorphism.
\end{proof}

The proposition above basically completes the proof of idempotency of $\mc{D}_{h,B}$, as we shall see below momentarily. But first, we establish this Proposition in the generality of twisted distributions. To do this, we shall use filtration and grading techniques, for which we shall first need the following lemma.

\begin{lemma}
There exists an exact sequence of $\Bplusbs[G_0]'$-modules of the form
$$\prod_J \Bplusbs[G_0]' \to \prod_I \Bplusbs[G_0]' \to \mc{D}_{h,B^+}' \to 0.$$
\end{lemma}

\begin{proof}
Take $I = \ZZ_{\ge 0}^d$ and $J = I \times I$. Let $1_t$ denote the element which is $1$ in the entry indexed by $t$ and $0$ elsewhere. The map $\prod_I \Bplusbs[G_0]' \to \mc{D}_{h,B^+}'$ 
is given by mapping $1_{\ul{n}} \mapsto \varpi^{-v_h(\ul{n})}{\bf{c}}^{\ul{n}}$, and is clearly surjective. The relations of the 
form ${\bf{c}}^{\ul{m}}\cdot \varpi^{-v_h(\ul{n})}{\bf{c}}^{\ul{n}} = \sum_{\ul{k}}a_{\ul{k}}(\ul{n},\ul{m})\varpi^{-v_h(\ul{k})}{\bf{c}}^{\ul{k}}$. Thus mapping $1_{\ul{n},\ul{m}}$ to ${\bf{c}}^{\ul{m}}\cdot 1_{\ul{n}} - (a_{\ul{k}}(\ul{n},\ul{m}))_{\ul{k}}$ gives the exact sequence.
\end{proof}

\begin{corollary}
The map $$H_0(\varpi\mc{D}_{h,B^+}'\otimes^{L}_{\Bplusbs[G_0]}\mc{D}_{h,B^+}') \to \varpi \mc{D}_{h,B^+}'$$ is an isomorphism.
\end{corollary}

\begin{proof}
By the previous lemma and Lemma \ref{lem:tens_prod_of_prod_with_G0}, there are some index sets $I, J$ such that there is an exact sequence
$$\prod_J \Bplusbs[G_0]' \to \prod_I \Bplusbs[G_0]'\to H_0(\varpi\mc{D}_{h,B^+}'\otimes^{L}_{\Bplusbs[G_0]}\mc{D}_{h,B^+}') \to 0.$$
In particular, since $\Bplusbs[G_0]'$ is a $\varpi$-adically complete $\Bplusbs$-module, it follows from \cite[Lemma 091U (1)]{stacks-project} that $H_0(\varpi\mc{D}_{h,B^+}'\otimes^{L}_{\Bplusbs[G_0]}\mc{D}_{h,B^+}')$ is a $\varpi$-adically complete $\Bplusbs$-module. Since $(g-1)(\varpi^nB^+)\subset \varpi^{n+1}B^+$ for $g \in G_0, n \in \ZZ_{\ge 0}$, it follows that for the $\varpi$-adic filtration, we have
$$\gr_{\varpi}(H_0(\varpi\mc{D}_{h,B^+}'\otimes^{L}_{\Bplusbs[G_0]}\mc{D}_{h,B^+}')) = \gr_{\varpi}(H_0(\varpi\mc{D}_{h,B^+}\otimes^{L}_{\Bplusbs[G_0]}\mc{D}_{h,B^+})).$$
Hence, we know that the map
$$H_0(\varpi\mc{D}_{h,B^+}'\otimes^{L}_{\Bplusbs[G_0]}\mc{D}_{h,B^+}')\to \mc{D}_{h,B^+}'$$
becomes isomorphism after applying the $\gr_{\varpi}$ functor. But this functor is conservative when applied to complete objects, by \cite[Lemma 5.2 (i)]{bhatt2019topological}, so the map must be an isomorphic before applying the functor $\gr_{\varpi}$. This concludes the proof.
\end{proof}

Using this corollary we can conclude the proof of Theorem \ref{thm:idempotent_dists}. Indeed, tensoring with $\otimes_{\Bplusbs}B_\bs$ we get that the map 
$$H_0(\mc{D}_{h,B}'\otimes^{L}_{\Bbs[G_0]}\mc{D}_{h,B}')\to \mc{D}_{h,B}'$$
is an isomorphism, but since the map $\Bbs[G_0]' \to \mc{D}_{h,B}'$ is flat (Corollary \ref{cor:flatness_of_dists_over_measures}), the tensor product $\mc{D}_{h,B}'\otimes^{L}_{\Bbs[G_0]}\mc{D}_{h,B}'$ is concentrated in degree 0.


\subsection{Characterization of locally analytic representations}

The work done in the previous subsection allow us to show that locally analytic representations sits nicely inside the category of continuous representations. We follow the treatment of \cite[{$\mathsection{4.3}$}]{RJRC22}.
\begin{theorem}
\label{thm:full_faithfullness}
The categories $\Mod^{\mt{solid}}_{\mc{D}_{h\han}(G,B)'}$, $\Mod^{\mt{solid}}_{\mc{D}^{h^+\han}(G,B)'}$ (respectively $D(\mc{D}_{h \han}(G,B)')$, $D(\mc{D}^{h^+\han}(G,B)')$) are full subcategories of $\Mod^{\mt{solid}}_{B_{\bs}[G]'}$ (respectively $D(B_{\bs}[G]')$).
\end{theorem}
\begin{proof}
This follows from idempotency (Theorem \ref{thm:idempotent_dists}). We give the proof for $\mc{D}_{h\han}(G,B)'$, the proof for $\mc{D}^{h^+\han}(G,B)'$ is similar. Namely, given $C,C'\in D(\mc{D}_{h\han}(G,B)')$, we may compute:
\begin{align*}
 &\ul{\mt{RHom}}_{B_{\bs}[G]'}(C,C')\\
 &=\ul{\mt{RHom}}_{\mc{D}_{h\han}(G,B)'}(\mc{D}_{h\han}(G,B)'\otimes^L_{B_{\bs}[G]'} C,C') \\&=\ul{\mt{RHom}}_{\mc{D}_{h\han}(G,B)'}(\mc{D}_{h\han}(G,B)'\otimes^L_{B_{\bs}[G]'}(\mc{D}_{h\han}(G,B)'\otimes^L_{\mc{D}_{h\han}(G,B)'} C),C')\\
&=\ul{\mt{RHom}}_{\mc{D}_{h\han}(G,B)'}((\mc{D}_{h\han}(G,B)'\otimes^L_{B_{\bs}[G]'}\mc{D}_{h\han}(G,B)')\otimes^L_{\mc{D}_{h\han}(G,B)'} C,C')\\
&=\ul{\mt{RHom}}_{\mc{D}_{h\han}(G,B)'}(\mc{D}_{h\han}(G,B)'\otimes^L_{\mc{D}_{h\han}(G,B)'} C,C')\\
&=\ul{\mt{RHom}}_{\mc{D}_{h\han}(G,B)'}(C,C'),
\end{align*}
which concludes the proof. \end{proof}


The following theorem then shows that locally analytic representations are in a sense the same as modules over distribution algebras.

\begin{theorem}
\label{thm:analytic_reps_characterization}
1. The category of lower $h$-analytic (respectively $h^+$-analytic) semilinear $G$-representations over $B$ is equal to the category of solid modules over $\mc{D}_{h \han}(G_0,B)'$ (respectively over $\mc{D}^{h^+ \han}(G_0,B)'$).

2. A complex $C \in D(B_{\bs}[G]')$ is lower $h$-analytic (respectively $h^+$-analytic) if and only if for all $n \in \ZZ$ the cohomologies $H^n(C)$ are lower $h$-analytic (respectively $h^+$-analytic). Equivalently, $C$ is in the essential image of $D(\mc{D}_{h \han}(G_0,B)')$ (respectively of $D(\mc{D}^{h^+ \han}(G_0,B)')$).
\end{theorem}

\begin{proof}
We shall provide proofs for $\mc{D}_{h \han}(G,B)'$. As usual the proofs for $\mc{D}^{h^+ \han}(G,B)'$ are similar.

 1. By definition, a lower $h$-analytic representation $V$ satisfies $$\ul{\mt{Hom}}_{B_{\bs}[G]'}(\mc{D}_{h \han}(G,B)',V)=V.$$
Acting with $\mc{D}_{h \han}(G,B)'$ on itself by right multiplication gives therefore gives $V$ the structure of a $\mc{D}_{h \han}(G,B)'$-module. Conversely, suppose $V$ is a $\mc{D}_{h \han}(G,B)'$-module. We have 
\begin{align*}
\ul{\mt{Hom}}_{B_{\bs}[G]'}(\mc{D}_{h \han}(G,B)',V)&=\ul{\mt{Hom}}_{\mc{D}_{h \han}(G,B)'}(\mc{D}_{h \han}(G,B)'\otimes^L_{B_{\bs}[G]'}\mc{D}_{h \han}(G,B)',V)\\
&=\ul{\mt{Hom}}_{\mc{D}_{h \han}(G,B)'}(\mc{D}_{h \han}(G,B)',V)=V,
\end{align*}
as required. Here for the second equality we used Theorem \ref{thm:idempotent_dists}.

2. We can argue as in part 1, replacing $\ul{\mt{Hom}}$ with $\ul{\mt{RHom}}$. For the cohomological characterization, use Theorem \ref{thm:solid_modules}.2. 
\end{proof}


\subsection{Comparison of cohomology.}

\begin{theorem}
\label{thm:cohomological_comparison}
Let $C\in D(B_{\bs}[G]')$ be a complex, then:

1. For every $h > 0,$ we have
$$\ul{\mt{RHom}}_{B_{\bs}[G]'}(B,C)=\ul{\mt{RHom}}_{B_{\bs}[G]'}(B,C_{h \han}).$$

2. For every $h > 0,$ we have
$$\ul{\mt{RHom}}_{B_{\bs}[G]'}(B,C)=\ul{\mt{RHom}}_{B_{\bs}[G]'}(B,C^{h^+ \han}).$$

3. We have
$$\ul{\mt{RHom}}_{B_{\bs}[G]'}(B,C)=\ul{\mt{RHom}}_{B_{\bs}[G]'}(B,C^{\la}).$$
\end{theorem}
\begin{proof}
Since $B$ is compact object as a $B_{\bs}[G]'$-module (as follows from Theorem \ref{thm:lazard_serre_resolution}), part 3 follows from either part 1 or 2. We prove part 1, part 2 is proved similarly. We compute
\begin{align*}
\ul{\mt{RHom}}_{B_{\bs}[G]'}(B,C_{h \han}) &= \ul{\mt{RHom}}_{B_{\bs}[G]'}(B,\ul{\mt{RHom}}_{B_{\bs}[G]'}(\mc{D}_{h\han}(G,B)',C))\\
&=\ul{\mt{RHom}}_{B_{\bs}[G]'}(\mc{D}_{h \han}(G,B)'\otimes_{B_{\bs}[G]'}^LB,C)\\
&=\ul{\mt{RHom}}_{B_{\bs}[G]'}(B,C),
\end{align*}
where in the last equality we used Corollary \ref{cor:dist_tens_triv}.
\end{proof}

In particular, we prove a (generalization of) conjectures 3.4 and 3.5 of \cite{Po24}. 

\begin{corollary}
Let $V$ be a $B$-Banach space. Then for the derived locally analytic vectors $\mt{R}^i_{\la}(V)$, there exists a spectral sequence 
$$E_2^{i,j} = \mt{Ext}_{B_{\bs}[G]'}^i(B,\mt{R}^j_{\la}(V)) \Longrightarrow \mt{Ext}_{B_{\bs}[G]'}^{i+j}(B,V).$$
\end{corollary}

This has yet another corollary:
\begin{corollary}
Let $V$ be a $B$-Banach space with $\mt{R}^i_{\la}(V) = 0$ for $i\geq 1$. Then 
$$\mt{H}^i(G,V)=\mt{H}^i(G,V^{\la}).$$
\end{corollary}

\bibliography{main}

\begin{thebibliography}{DDSMS03}

\bibitem[AIP18]{andreatta2018halo}
Fabrizio Andreatta, Adrian Iovita, and Vincent Pilloni.
\newblock Le halo spectral.
\newblock In {\em Annales scientifiques de l'{\'E}cole Normale Sup{\'e}rieure}, volume~51, pages 603--655, 2018.

\bibitem[Ami64]{amice1964interpolation}
Yvette Amice.
\newblock Interpolation $ p $-adique.
\newblock {\em Bulletin de la Soci{\'e}t{\'e} math{\'e}matique de France}, 92:117--180, 1964.

\bibitem[And21]{andreychev2021pseudocoherent}
Grigory Andreychev.
\newblock Pseudocoherent and perfect complexes and vector bundles on analytic adic spaces.
\newblock {\em arXiv preprint arXiv:2105.12591}, 2021.

\bibitem[BC16]{berger2016theorie}
Laurent Berger and Pierre Colmez.
\newblock Th{\'e}orie de {Sen} et vecteurs localement analytiques.
\newblock {\em Ann. Sci. {\'E}c. Norm. Sup{\'e}r.(4)}, 49(4):947--970, 2016.

\bibitem[Bel24a]{bellovin2024cohomology}
Rebecca Bellovin.
\newblock Cohomology of ({$\varphi$}, {$\Gamma$})-modules over pseudorigid spaces.
\newblock {\em International Mathematics Research Notices}, 2024(4):2999--3051, 2024.

\bibitem[Bel24b]{bellovin2024modularity}
Rebecca Bellovin.
\newblock Modularity of trianguline galois representations.
\newblock In {\em Forum of Mathematics, Sigma}, volume~12, page~e3. Cambridge University Press, 2024.

\bibitem[Ber02]{berger2002representations}
Laurent Berger.
\newblock Repr{\'e}sentations {$p$}-adiques et {\'e}quations diff{\'e}rentielles.
\newblock {\em Inventiones mathematicae}, 148(2):219--284, 2002.

\bibitem[Ber16]{berger2016multivariable}
Laurent Berger.
\newblock Multivariable {($\varphi, \Gamma$)}-modules and locally analytic vectors.
\newblock {\em Duke Mathematical Journal}, 165(18):3567, 2016.

\bibitem[BMS19]{bhatt2019topological}
Bhargav Bhatt, Matthew Morrow, and Peter Scholze.
\newblock Topological {Hochschild} homology and integral {$p$}-adic {Hodge} theory.
\newblock {\em Publications math{\'e}matiques de l'IH{\'E}S}, 129(1):199--310, 2019.

\bibitem[BR22]{BR22a}
Laurent Berger and Sandra Rozensztajn.
\newblock Decompletion of cyclotomic perfectoid fields in positive characteristic.
\newblock {\em Annales Henri Lebesgue}, 5:1261--1276, 2022.

\bibitem[BR24]{BR22b}
Laurent Berger and Sandra Rozensztajn.
\newblock {Super}-{H{\"o}lder} vectors and the field of norms.
\newblock {\em Algebra \& Number Theory}, 19(1):195--211, 2024.

\bibitem[Cam22]{camargo2022geometric}
JE~Camargo.
\newblock Geometric {Sen} theory over rigid analytic spaces.
\newblock {\em arXiv preprint arXiv:2205.02016}, 2022.

\bibitem[CC98]{cherbonnier1998representations}
Fr{\'e}d{\'e}ric Cherbonnier and Pierre Colmez.
\newblock Repr{\'e}sentations {$p$}-adiques surconvergentes.
\newblock {\em Inventiones mathematicae}, 133(3):581--611, 1998.

\bibitem[CM98]{coleman1998eigencurve}
Robert Coleman and Barry Mazur.
\newblock The eigencurve.
\newblock {\em London Mathematical Society Lecture Note Series}, pages 1--114, 1998.

\bibitem[Col10]{colmez2010representations}
Pierre Colmez.
\newblock Repr{\'e}sentations de {$\mathrm{GL}_2(\QQ_p)$} et ({$\varphi$}, {$\Gamma$})-modules.
\newblock {\em Ast{\'e}risque}, 330(281):509, 2010.

\bibitem[Col16]{colmez2016representations}
Pierre Colmez.
\newblock Repr{\'e}sentations localement analytiques de {$\mathrm{GL}_2(\QQ_p)$} et ({$\varphi$}, {$\Gamma$})-modules.
\newblock {\em Representation Theory of the American Mathematical Society}, 20(9):187--248, 2016.

\bibitem[CS19a]{CS19an}
Dustin Clausen and Peter Scholze.
\newblock Analytic geometry.
\newblock {\em Lecture notes. Available at Scholze’s webpage}, 2019.

\bibitem[CS19b]{CS19con}
Dustin Clausen and Peter Scholze.
\newblock Condensed mathematics.
\newblock {\em Lecture notes. Available at Scholze’s webpage}, 2019.

\bibitem[DDSMS03]{DdMS03}
John~D Dixon, Marcus~PF Du~Sautoy, Avinoam Mann, and Dan Segal.
\newblock {\em Analytic pro-{$p$} groups}.
\newblock Number~61. Cambridge University Press, 2003.

\bibitem[Dos11]{dospinescu2011equations}
Gabriel Dospinescu.
\newblock Equations diff{\'e}rentielles {$p$}-adiques et modules de {Jacquet} analytiques.
\newblock {\em Automorphic Forms and Galois Representations}, 1:359--374, 2011.

\bibitem[EGH22]{emerton2022introduction}
Matthew Emerton, Toby Gee, and Eugen Hellmann.
\newblock An introduction to the categorical {$p$}-adic langlands program.
\newblock {\em arXiv preprint arXiv:2210.01404}, 2022.

\bibitem[Eme06]{emerton2006interpolation}
Matthew Emerton.
\newblock On the interpolation of systems of eigenvalues attached to automorphic hecke eigenforms.
\newblock {\em Inventiones mathematicae}, 164(1):1--84, 2006.

\bibitem[Eme11]{emerton2011local}
Matthew Emerton.
\newblock Local-global compatibility in the $p$-adic langlands programme for {$\mathrm{GL}_{2_{/\QQ}}$}.
\newblock {\em preprint}, 3(4), 2011.

\bibitem[Eme17]{emerton2017locally}
Matthew Emerton.
\newblock {\em Locally analytic vectors in representations of locally {$p$}-adic analytic groups}, volume 248.
\newblock American Mathematical Society, 2017.

\bibitem[GP21]{gao2021locally}
Hui Gao and L{\'e}o Poyeton.
\newblock Locally analytic vectors and overconvergent ({$\varphi$},{$\tau$})-modules.
\newblock {\em Journal of the Institute of Mathematics of Jussieu}, 20(1):137--185, 2021.

\bibitem[Gul19]{gulotta2019equidimensional}
Daniel~R Gulotta.
\newblock Equidimensional adic eigenvarieties for groups with discrete series.
\newblock {\em Algebra \& Number Theory}, 13(8):1907--1940, 2019.

\bibitem[HN17]{hansen2017universal}
David Hansen and James Newton.
\newblock Universal eigenvarieties, trianguline {Galois} representations, and {$p$}-adic {Langlands} functoriality.
\newblock {\em Journal f{\"u}r die reine und angewandte Mathematik}, 2017(730):1--64, 2017.

\bibitem[HVO96]{huishi1996zariskian}
Li~Huishi and Freddy Van~Oystaeyen.
\newblock Zariskian filtrations.
\newblock In {\em Zariskian Filtrations}, pages 70--126. Springer, 1996.

\bibitem[JC23]{RJRC23}
Joaqu{\'\i}n~Rodrigues Jacinto and Juan Esteban~Rodr{\'\i}guez Camargo.
\newblock Solid locally analytic representations.
\newblock {\em arXiv preprint arXiv:2305.03162}, 2023.

\bibitem[JN19]{johansson2019extended}
Christian Johansson and James Newton.
\newblock Extended eigenvarieties for overconvergent cohomology.
\newblock {\em Algebra \& Number Theory}, 13(1):93--158, 2019.

\bibitem[Ked04]{kedlaya2004p}
Kiran~S Kedlaya.
\newblock A {$p$}-adic local monodromy theorem.
\newblock {\em Annals of mathematics}, pages 93--184, 2004.

\bibitem[Kis03]{kisin2003overconvergent}
Mark Kisin.
\newblock Overconvergent modular forms and the {Fontaine}-{Mazur} conjecture.
\newblock {\em Inventiones mathematicae}, 153(2):373--454, 2003.

\bibitem[Koh11]{kohlhaase2011cohomology}
Jan Kohlhaase.
\newblock The cohomology of locally analytic representations.
\newblock 2011.

\bibitem[Laz65]{lazard1965groupes}
Michel Lazard.
\newblock Groupes analytiques $ p $-adiques.
\newblock {\em Publications Math{\'e}matiques de l'IH{\'E}S}, 26:5--219, 1965.

\bibitem[Lou17]{Lou17}
Joao~NP Louren{\c{c}}o.
\newblock The {Riemannian} {Hebbarkeitssätze} for pseudorigid spaces.
\newblock {\em arXiv preprint arXiv:1711.06903}, 2017.

\bibitem[Man22]{Ma22}
Lucas Mann.
\newblock {\em A $p$-adic 6-functor formalism in rigid-analytic geometry}.
\newblock PhD thesis, Rheinische Friedrich-Wilhelms-Universitaet Bonn (Germany), 2022.

\bibitem[Pan22]{pan2022locally}
Lue Pan.
\newblock On locally analytic vectors of the completed cohomology of modular curves.
\newblock In {\em Forum of Mathematics, Pi}, volume~10, page~e7. Cambridge University Press, 2022.

\bibitem[Por24a]{porat2024locally}
Gal Porat.
\newblock Locally analytic vector bundles on the fargues--fontaine curve.
\newblock {\em Algebra \& Number Theory}, 18(5):899--946, 2024.

\bibitem[Por24b]{Po24}
Gal Porat.
\newblock Locally analytic vectors and decompletion in mixed characteristic.
\newblock {\em arXiv preprint arXiv:2407.19791}, 2024.

\bibitem[RJRC22]{RJRC22}
Joaqu{\'\i}n Rodrigues~Jacinto and Juan Rodr{\'\i}guez~Camargo.
\newblock Solid locally analytic representations of $p$-adic lie groups.
\newblock {\em Representation Theory of the American Mathematical Society}, 26(31):962--1024, 2022.

\bibitem[Ser09]{serre2009lie}
Jean-Pierre Serre.
\newblock {\em Lie algebras and Lie groups: 1964 lectures given at Harvard University}.
\newblock Springer, 2009.

\bibitem[ST01]{ST01}
Peter Schneider and Jeremy Teitelbaum.
\newblock {$p$}-adic {Fourier} theory.
\newblock {\em Documenta Mathematica}, 6:447--481, 2001.

\bibitem[ST02a]{schneider_teitelbaum2002algebras}
Peter Schneider and Jeremy Teitelbaum.
\newblock Algebras of $p$-adic distributions and admissible representations.
\newblock {\em arXiv preprint math/0206056}, 2002.

\bibitem[ST02b]{ST02}
Peter Schneider and Jeremy Teitelbaum.
\newblock Algebras of $p$-adic distributions and admissible representations.
\newblock {\em arXiv preprint math/0206056}, 2002.

\bibitem[ST02c]{ST02locally}
Peter Schneider and Jeremy Teitelbaum.
\newblock Locally analytic distributions and {$p$}-adic representation theory, with applications to {$\mathrm{GL}_2$}.
\newblock {\em Journal of the American Mathematical Society}, 15(2):443--468, 2002.

\bibitem[{Sta}25]{stacks-project}
The {Stacks project authors}.
\newblock The stacks project.
\newblock \url{https://stacks.math.columbia.edu}, 2025.

\end{thebibliography}
\end{document}